# EXISTENCE OF INDEPENDENT RANDOM MATCHING

By Darrell Duffie and Yeneng Sun

*Stanford University and National University of Singapore*

This paper shows the existence of independent random matching of a large (continuum) population in both static and dynamic systems, which has been popular in the economics and genetics literatures. We construct a joint agent-probability space, and randomized mutation, partial matching and match-induced type-changing functions that satisfy appropriate independence conditions. The proofs are achieved via nonstandard analysis. The proof for the dynamic setting relies on a new Fubini-type theorem for an infinite product of Loeb transition probabilities, based on which a continuum of independent Markov chains is derived from random mutation, random partial matching and random type changing.

**1. Introduction.** Economists and geneticists, among others, have implicitly or explicitly assumed the exact law of large numbers for independent random matching in a continuum population, by which we mean a nonatomic measure space of agents. This result is relied upon in large literatures within general equilibrium theory, game theory, monetary theory, labor economics, illiquid financial markets and biology, as discussed in [15], which provides extensive references. Such a law of large numbers allows a dramatic simplification of the dynamics of the cross-sectional distribution of properties among a large population of randomly interacting agents. Mathematical foundations, however, have been lacking, as has been noted by Green and Zhou [19].

Given the fundamental measurability problems associated with modeling a continuum of independent random variables,[1] there has, up to now, been no theoretical treatment of the exact law of large numbers for independent random matching among a continuum population. In [40], various versions of the exact law of large numbers and their converses are proved by applying simple measure-theoretic methods to an extension of the usual product probability space that has the Fubini property.[2] The measure-theoretic framework of [40] is adopted in our companion paper [15] to obtain an exact law of large numbers for random pairwise matching by formulating a suitable independence condition in types.[3]

[1]See, for example, [12, 13, 24] and detailed discussions in [38] and [40].

In particular, assuming independent random matching in a continuum population, and some related independence assumptions regarding random mutation and match-induced type changing, we prove in [15] that there is an almost sure deterministic cross-sectional distribution of types in a large population for both static and dynamic systems, a property that had been widely used without a proper foundation. In addition, we show in [15] that the time evolution of the cross-sectional distribution of types can be completely determined from the agent-level Markov chain for type, with known transition matrices.

The main aim of this paper is to provide the first theoretical treatment of the existence of independent random matching in a continuum population. In particular, we construct a joint agent-probability space, and randomized mutation, partial matching and match-induced type-changing functions that satisfy the independence conditions in [15]. Though the existence results of this paper are stated using common measure-theoretic terms, their proofs make extensive use of nonstandard analysis. One can pick up some background knowledge on nonstandard analysis from the first three chapters of the book [34].

Since our dynamic system of independent random matching generates a dependence structure across time (in particular, a continuum of independent Markov chains of agents' type), we need to construct an infinite product of Loeb transition probabilities. Our proof of the existence result for the dynamic setting is based on such an infinite product and an associated Fubini-type theorem that is derived from Keisler [26]. Specifically, in order to prove the Fubini-type property of an infinite product of Loeb transition probabilities, we first generalize Keisler's Fubini theorem for the product of two Loeb probability measures, from [26], to the setting of a Loeb transition probability (i.e., a class of Loeb measures) in Theorem 5.1. In constructing the internal transition probabilities for the step of random mutation in Section 6, we use a hyperfinite product space and its coordinate functions, which are closely related to those of Keisler [26], and in particular to the law of large numbers for a hyperfinite sequence of $*$-independent random variables, as noted in [26], page 56.

---

[2]These results were originally stated on Loeb measure spaces in [38]. However, as noted in [40], they can be proved for an extension of the usual product space that has the Fubini property; see also Chapter 7 in [34] (and in particular, Sections 7.5 and 7.6), written by Sun.

[3]The independence condition that we propose in [15] is natural, but may not be obvious. For example, random matching in a finite population may not allow independence among agents since the matching of agent $i$ to agent $j$ implies of course that $j$ is also matched to $i$, implying some correlation among agents. The effect of this correlation is reduced to zero in a continuum population. In other words, one may obtain standard independence by rounding infinitesimals but not ∗-independence for a hyperfinite population (see Section 4 below). A new concept, "Markov conditional independence in types," is proposed in [15] for dynamic matching, under which the transition law at each randomization step depends on only the previous one or two steps of randomization.

Historically, reliance on the exact law of large numbers for independent random matching dates back at least to 1908, when G. H. Hardy [20] and W. Weinberg (see [6]) independently proposed that random mating in a large population leads to constant and easily calculated fractions of each allele in the population. Hardy wrote: "suppose that the numbers are fairly large, so that the mating may be regarded as random," and then used, in effect, an exact law of large numbers for random matching to deduce his results. Consider, for illustration, a continuum population of gametes consisting of two alleles, $A$ and $B$, in initial proportions $p$ and $q = 1 - p$. Then, following the Hardy–Weinberg approach, the new population would have a fraction $p^2$ whose parents are both of type $A$, a fraction $q^2$ whose parents are both of type $B$, and a fraction $2pq$ whose parents are of mixed type (heterozygotes). These genotypic proportions asserted by Hardy and Weinberg are already, implicitly, based on an exact law of large numbers for independent random matching in a large population. We provide a suitable existence framework.

Going from a static to a dynamic environment, we also provide an existence result that allows the computation of a steady-state constant deterministic population distribution of types. For illustration, suppose in the Hardy–Weinberg setting above that with both parents of allele $A$, the offspring are of allele $A$, and that with both parents of allele $B$, the offspring are of allele $B$. Suppose that the offspring of parents of different alleles are, say, equally likely to be of allele $A$ or allele $B$. The Hardy–Weinberg equilibrium for this special case is a population with steady-state constant proportions $p = 60\%$ of allele $A$ and $q = 40\%$ of allele $B$. With the law of large numbers for independent random matching, this is verified by checking that, if generation $k$ has this cross-sectional distribution, then the fraction of allele $A$ in generation $k + 1$ is almost surely $0.6^2 + 0.5 \times (2 \times 0.6 \times 0.4) = 0.6$. Our existence results for a dynamic model of random matching provide a mathematical foundation for this Hardy–Weinberg law governing steady-state allelic and genotypic frequencies.

In the field of economics, Hellwig [21] is the first, to our knowledge, to have relied on the effect of the exact law of large numbers for random pairwise matching in a market, in a 1976 study of a monetary exchange economy.[4] Since the 1970s, a large economics literature has routinely relied on an exact law of large numbers for independent random matching in a continuum population. This implicit use of this result occurs in general equilibrium theory (e.g., [17, 18, 35, 41]), game theory (e.g., [4, 5, 8, 16, 22]), monetary theory (e.g., [11, 19, 21, 31]), labor economics (e.g., [10, 23, 36, 37]) and financial market theory (e.g., [14, 32]). Mathematical foundations, however, have been lacking, as has been noted by Green and Zhou [19]. In almost all of this literature, dynamics are crucial. For example, in the monetary and finance literature cited above, each agent in the economy solves a dynamic programming problem that is based in part on the conjectured dynamics of the cross-sectional distribution of agent types. An equilibrium has the property that the combined effect of individually optimal dynamic behavior is consistent with the conjectured population dynamics. In order to simplify the analysis, much of the literature relies on equilibria with a stationary distribution of agent types.

The remainder of the paper is organized as follows. Section 2 considers existence of independent random matching, both full and partial, in a static setting, after a brief introduction of the measure-theoretic framework in Section 2.1. Theorem 2.4 of Section 2.2 shows the existence of random full matching with independent types, meaning roughly that, for essentially every pair $(i, j)$ of agents, the type of the agent to be randomly matched with agent $i$ is independent of the type of the agent to be randomly matched with agent $j$. Theorem 2.6 of Section 2.3 then considers existence for the case of random search and matching, that is, for random partial matchings that are independent in types. Proofs of Theorems 2.4 and 2.6, which use nonstandard analysis extensively in the computations, are given in Section 4.

Section 3 considers a dynamical system for agent types, allowing for random mutation, partial matching and match-induced random type changes. We borrow from our companion paper [15] the inductive definition of such a dynamical system given in Section 3.1, and the condition of Markov conditional independence found in Section 3.2. The latter condition captures the idea that at every time period, there are three stages: (1) an independent random mutation, (2) an independent random partial matching, and

---

[4]Diamond [9] had earlier treated random matching of a large population with, in effect, finitely many employers, but not pairwise matching within a large population. The matching of a large population with a finite population can be treated directly by the exact law of large numbers for a continuum of independent random variables. For example, let $N(i)$ be the event that worker $i$ is matched with an employer of a given type, and suppose this event is pairwise independent and of the same probability $p$, in a continuum population of such workers. Then, under the conditions of [40], the fraction of the population that is matched to this type of employer is $p$, almost surely.

(3) for those agents matched, an independent random type change induced by matching. In economics applications, for example, match-induced type changes arise from productivity shocks, changes in asset holdings induced by trade between the matched agents, changes in credit positions, or changes in money holdings. Theorem 3.1 of Section 3.3 shows the existence of a dynamical system $\mathbb{D}$ with random mutation, partial matching and type changing that is Markov conditionally independent in types with any given parameters. Theorem 3 of [15] then implies that the type processes of individual agents in such a dynamical system $\mathbb{D}$ form a continuum of independent Markov chains, and that the time evolution of the cross-sectional distribution of types is deterministic and completely determined from a Markov chain with explicitly calculated transition matrices.[5]

We prove Theorems 2.4 and 2.6 in Section 4. Turning to Section 5, we first prove in Section 5.1 a generalized Fubini theorem for a Loeb transition probability. Then, a generalized Ionescu–Tulcea theorem for an infinite sequence of Loeb transition probabilities is shown in Section 5.2. Finally, a Fubini extension based on Loeb product transition probability systems is constructed in Section 5.3. Based on the results in Section 5, we prove Theorem 3.1 in Section 6.

Finally, we emphasize again that we must work with extensions of the usual product measure spaces (of agents and states of the world), since a process formed by a continuum of independent random variables is never measurable with respect to the completion of the usual product $\sigma$-algebra, except in the trivial case that almost all of the random variables in the process are constants.[6]

**2. Existence of independent random matchings in the static case.** In this section, we first give some background definitions in Section 2.1. Then, we consider the existence of random matchings that are independent in types, for full and partial matchings, in Sections 2.2 and 2.3, respectively. Proofs of the two existence theorems, Theorems 2.4 and 2.6, are given in Sections 4.1 and 4.2. The full and partial matching models in Theorems 2.4 and 2.6 satisfy the respective conditions in Theorems 1 and 2 in [15], which means that the respective conclusions of Theorems 1 and 2 in [15], characterizing the implications of the law of large numbers for process describing the cross-sectional distribution of types, hold for these matching models.

---

[5]The models in [27, 28] and [29] assume that there is a set of individual agents and a single marketmaker. At each time, a randomly chosen individual agent trades with the marketmaker. References [1] and [2] formalize a link between matching and informational constraints, which do not consider random matching under the independence assumption as in the model here.

[6]See, for example, Proposition 1.1 in [39].

2.1. *Some background definitions.* Let probability spaces $(I, \mathcal{I}, \lambda)$ and $(\Omega, \mathcal{F}, P)$ be our index and sample spaces, respectively.[7] In our applications, $(I, \mathcal{I}, \lambda)$ is an atomless probability space that indexes the agents.[8] Let $(I \times \Omega, \mathcal{I} \otimes \mathcal{F}, \lambda \otimes P)$ be the usual product probability space. For a function $f$ on $I \times \Omega$ (not necessarily $\mathcal{I} \otimes \mathcal{F}$-measurable), and for $(i, \omega) \in I \times \Omega$, $f_i$ represents the function $f(i, \cdot)$ on $\Omega$, and $f_\omega$ the function $f(\cdot, \omega)$ on $I$.

In order to work with independent type processes arising from random matching, we need to work with an extension of the usual measure-theoretic product that retains the Fubini property. A formal definition, as in [40], is as follows.

DEFINITION 2.1. A probability space $(I \times \Omega, \mathcal{W}, Q)$ extending the usual product space $(I \times \Omega, \mathcal{I} \otimes \mathcal{F}, \lambda \otimes P)$ is said to be a *Fubini extension* of $(I \times \Omega, \mathcal{I} \otimes \mathcal{F}, \lambda \otimes P)$ if for any real-valued $Q$-integrable function $g$ on $(I \times \Omega, \mathcal{W})$, the functions $g_i = g(i, \cdot)$ and $g_\omega = f(\cdot, \omega)$ are integrable respectively on $(\Omega, \mathcal{F}, P)$ for $\lambda$-almost all $i \in i$ and on $(I, \mathcal{I}, \lambda)$ for $P$-almost all $\omega \in \Omega$; and if, moreover, $\int_\Omega g_i \, dP$ and $\int_I g_\omega \, d\lambda$ are integrable respectively on $(I, \mathcal{I}, \lambda)$ and on $(\Omega, \mathcal{F}, P)$, with $\int_{I \times \Omega} g \, dQ = \int_I (\int_\Omega g_i \, dP) \, d\lambda = \int_\Omega (\int_I g_\omega \, d\lambda) \, dP$. To reflect the fact that the probability space $(I \times \Omega, \mathcal{W}, Q)$ has $(I, \mathcal{I}, \lambda)$ and $(\Omega, \mathcal{F}, P)$ as its marginal spaces, as required by the Fubini property, it will be denoted by $(I \times \Omega, \mathcal{I} \boxtimes \mathcal{F}, \lambda \boxtimes P)$.

An $\mathcal{I} \boxtimes \mathcal{F}$-measurable function $f$ will also be called a process, while $f_i$ is called a random variable of the process and $f_\omega$ is called a sample function of the process.

We now introduce the following crucial independence condition. We state the definition using a complete separable metric space $X$ for the sake of generality; in particular, a finite space or a Euclidean space is a complete separable metric space.

DEFINITION 2.2. An $\mathcal{I} \boxtimes \mathcal{F}$-measurable process $f$ from $I \times \Omega$ to a complete separable metric space $X$ is said to be essentially pairwise independent if for $\lambda$-almost all $i \in I$, the random variables $f_i$ and $f_j$ are independent for $\lambda$-almost all $j \in I$.

---

[7]All measures in this paper are countably additive set functions defined on $\sigma$-algebras.

[8]A probability space $(I, \mathcal{I}, \lambda)$ is atomless if there does not exist $A \in \mathcal{I}$ such that $\lambda(A) > 0$, and for any $\mathcal{I}$-measurable subset $C$ of $A$, $\lambda(C) = 0$ or $\lambda(C) = \lambda(A)$. For those interested in the case of a literal continuum of agents, it is noted in the beginning of Section 4 that one can indeed take $I$ to be the unit interval with some atomless probability measure. Corollary 4.3 in [40] shows, however, that, in general, it makes no sense to impose the Lebesgue measure structure when an independent process is considered.

2.2. *Existence of independent random full matchings.* We follow the notation in Section 2.1. Below is a formal definition of random full matching.

DEFINITION 2.3 (Full matching).

1. Let $S = \{1, 2, \ldots, K\}$ be a finite set of types, $\alpha : I \to S$ an $\mathcal{I}$-measurable type function of agents and $p$ its distribution on $S$. For $1 \leq k \leq K$, let $I_k = \{i \in I : \alpha(i) = k\}$ and $p_k = \lambda(I_k)$ for each $1 \leq k \leq K$.
2. A full matching $\phi$ is a bijection from $I$ to $I$ such that for each $i \in I$, $\phi(i) \neq i$ and $\phi(\phi(i)) = i$.
3. A random full matching $\pi$ is a mapping from $I \times \Omega$ to $I$ such that (i) $\pi_\omega$ is a full matching for each $\omega \in \Omega$; (ii) if we let $g$ be the type process $\alpha(\pi)$, then $g$ is measurable from $(I \times \Omega, \mathcal{I} \boxtimes \mathcal{F}, \lambda \boxtimes P)$ to $S$; (iii) for $\lambda$-almost all $i \in I$, $g_i$ has distribution $p$.
4. A random full matching $\pi$ is said to be independent in types if the type process $g$ is essentially pairwise independent.

Condition (1) of this definition says that a fraction $p_k$ of the population is of type $k$. Condition (2) says that there is no self-matching, and that if $i$ is matched to $j = \phi(i)$, then $j$ is matched to $i$. Condition (3)(iii) means that for almost every agent $i$, the probability that $i$ is matched to a type-$k$ agent is $p_k$, the fraction of type-$k$ agents in the population. Condition (4) says that for almost all agents $i$ and $j \in I$, the event that agent $i$ is matched to a type-$k$ agent is independent of the event that agent $j$ is matched to a type-$l$ agent for any $k, l \in S$.

The following theorem shows the existence of an independent random full matching model that satisfies a few strong conditions that are specified in footnote 4 of [35], and is universal in the sense that it does not depend on particular type functions.[9] Note that condition (1)(ii) below implies that for any $i, j \in I$, $P(\pi_i = j) = 0$ since $\lambda(\{j\}) = 0$, which means that the probability that agent $i$ is matched with a given agent $j$ is zero.

THEOREM 2.4. *There exists an atomless probability space $(I, \mathcal{I}, \lambda)$ of agents, a sample probability space $(\Omega, \mathcal{F}, P)$, a Fubini extension $(I \times \Omega, \mathcal{I} \boxtimes \mathcal{F}, \lambda \boxtimes P)$ of the usual product probability space, and a random full matching $\pi$ from $(I \times \Omega, \mathcal{I} \boxtimes \mathcal{F}, \lambda \boxtimes P)$ to $I$ such that* 2.3. *The existence of independent random partial matchings.* We shall now consider the case of random partial matchings. The following is a formal definition.

---

[9]When $(I, \mathcal{I}, \lambda)$ is taken to be the unit interval with the Borel algebra and Lebesgue measure, property (1)(iii) of Theorem 2.4 can be restated as "for $P$-almost all $\omega \in \Omega$, $\lambda(A_1 \cap \pi_\omega^{-1}(A_2)) = \lambda(A_1)\lambda(A_2)$ holds for any $A_1, A_2 \in \mathcal{I}$" by using the fact that the countable collection of rational intervals in $[0, 1]$ generates the Borel algebra. Footnote 4 of [35] shows the *nonexistence* of a random full matching $\pi$ that satisfies (i)–(iii) of part (1).

1. (i) *for each $\omega \in \Omega$, $\lambda(\pi_\omega^{-1}(A)) = \lambda(A)$ for any $A \in \mathcal{I}$*, (ii) *for each $i \in I$, $P(\pi_i^{-1}(A)) = \lambda(A)$ for any $A \in \mathcal{I}$*, (iii) *for any $A_1, A_2 \in \mathcal{I}$, $\lambda(A_1 \cap \pi_\omega^{-1}(A_2)) = \lambda(A_1)\lambda(A_2)$ holds for $P$-almost all $\omega \in \Omega$*;
2. *$\pi$ is independent in types with respect to any given type function $\alpha$ from $I$ to any finite type space $S$.*

DEFINITION 2.5. Let $\alpha: I \to S$ be an $\mathcal{I}$-measurable type function with type distribution $p = (p_1, \ldots, p_K)$ on $S$. Let $\pi$ be a mapping from $I \times \Omega$ to $I \cup \{J\}$, where $J$ denotes "no match."

1. We say that $\pi$ is a random partial matching with no-match probabilities $q_1, \ldots, q_K$ in $[0,1]$ if:

    (i) For each $\omega \in \Omega$, the restriction of $\pi_\omega$ to $I - \pi_\omega^{-1}(\{J\})$ is a full matching on $I - \pi_\omega^{-1}(\{J\})$.[10]

    (ii) After extending the type function $\alpha$ to $I \cup \{J\}$ such that $\alpha(J) = J$, and letting $g = \alpha(\pi)$, we have $g$ measurable from $(I \times \Omega, \mathcal{I} \boxtimes \mathcal{F}, \lambda \boxtimes P)$ to $S \cup \{J\}$.

    (iii) For $\lambda$-almost all $i \in I_k$, $P(g_i = J) = q_k$ and[11]
    $$P(g_i = l) = \frac{(1-q_k)p_l(1-q_l)}{\sum_{r=1}^K p_r(1-q_r)}.$$

2. A random partial matching $\pi$ is said to be independent in types if the process $g$ (taking values in $S \cup \{J\}$) is essentially pairwise independent.[12]

The following theorem generalizes Theorem 2.4 to the case of random partial matchings. Because the given parameters for no-matching probabilities may be type-dependent, it is not possible to produce a universal matching model for random partial matchings as in the case of full matchings.

---

[10] This means that an agent $i$ with $\pi_\omega(i) = J$ is not matched, while any agent in $I - \pi_\omega^{-1}(\{J\})$ is matched. This produces a partial matching on $I$.

[11] If an agent of type $k$ is matched, its probability of being matched to a type-$l$ agent should be proportional to the type distribution of matched agents. The fraction of the population of matched agents among the total population is $\sum_{r=1}^K p_r(1-q_r)$. Thus, the relative fraction of type-$l$ matched agents to that of all matched agents is $(p_l(1-q_l))/\sum_{r=1}^K p_r(1-q_r)$. This implies that the probability that a type-$k$ agent is matched to a type-$l$ agent is $(1-q_k)(p_l(1-q_l))/\sum_{r=1}^K p_r(1-q_r)$. When $\sum_{r=1}^K p_r(1-q_r) = 0$, we have $p_k(1-q_k) = 0$ for all $1 \leq k \leq K$, in which case almost no agents are matched, and we can interpret the ratio $((1-q_k)p_l(1-q_l))/\sum_{r=1}^K p_r(1-q_r)$ as zero.

[12] This means that for almost all agents $i, j \in I$, whether agent $i$ is unmatched or matched to a type-$k$ agent is independent of a similar event for agent $j$.

THEOREM 2.6. *There is an atomless probability space $(I, \mathcal{I}, \lambda)$ of agents such that for any given $\mathcal{I}$-measurable type function $\beta$ from $I$ to $S$, and for any $q \in [0,1]^S$, (1) there exists a sample space $(\Omega, \mathcal{F}, P)$ and a Fubini extension $(I \times \Omega, \mathcal{I} \boxtimes \mathcal{F}, \lambda \boxtimes P)$ of the usual product probability space; (2) there exists an independent-in-types random partial matching $\pi$ from $(I \times \Omega, \mathcal{I} \boxtimes \mathcal{F}, \lambda \boxtimes P)$ to $I$ with $q = (q_1, \ldots, q_K)$ as the no-match probabilities.*

**3. The existence of a dynamical system with random mutation, partial matching and type changing that is Markov conditionally independent in types.** A discrete-time dynamical system $\mathbb{D}$ with random mutation, partial matching and type changing that is Markov conditionally independent in types is introduced in Section 4 of our companion paper [15]. The purpose of this section is to show the existence of such a dynamical system $\mathbb{D}$ with any given parameters. In Sections 3.1 and 3.2, we reproduce respectively the inductive definition of a dynamical system with random mutation, partial matching and type changing and the condition of Markov conditional independence, which originated with [15]. The general existence of the dynamical system $\mathbb{D}$ is presented in Theorem 3.1 in Section 3.3 and its proof in Section 6.

3.1. *Definition of a dynamical system with random mutation, partial matching and type changing.* Let $S = \{1, 2, \ldots, K\}$ be a finite set of types. A discrete-time dynamical system $\mathbb{D}$ with random mutation, partial matching and type changing in each period can be defined intuitively as follows. The initial distribution of types is $p^0$. That is, $p^0(k)$ (denoted by $p_k^0$) is the initial fraction of agents of type $k$. In each time period, each agent of type $k$ first goes through a stage of random mutation, becoming an agent of type $l$ with probability $b_{kl}$. In models such as [14], for example, this mutation generates new motives for trade. Then, each agent of type $k$ is either not matched, with probability $q_k$, or is matched to a type-$l$ agent with a probability proportional to the fraction of type-$l$ agents in the population immediately after the random mutation step. When an agent is not matched, she keeps her type. Otherwise, when a pair of agents with respective types $k$ and $l$ are matched, each of the two agents changes types; the type-$k$ agent becomes type $r$ with probability $\nu_{kl}(r)$, where $\nu_{kl}$ is a probability distribution on $S$, and similarly for the type-$l$ agent.

We shall now define formally a dynamical system $\mathbb{D}$ with random mutation, partial matching and type changing. As in Section 2, let $(I, \mathcal{I}, \lambda)$ be an atomless probability space representing the space of agents, $(\Omega, \mathcal{F}, P)$ a sample probability space, and $(I \times \Omega, \mathcal{I} \boxtimes \mathcal{F}, \lambda \boxtimes P)$ a Fubini extension of the usual product probability space.

Let $\alpha^0 : I \to S = \{1, \ldots, K\}$ be an initial $\mathcal{I}$-measurable type function with distribution $p^0$ on $S$. For each time period $n \geq 1$, we first have a random mutation that is modeled by a process $h^n$ from $(I \times \Omega, \mathcal{I} \boxtimes \mathcal{F}, \lambda \boxtimes P)$ to $S$, then a random partial matching described by a function $\pi^n$ from $(I \times \Omega, \mathcal{I} \boxtimes \mathcal{F}, \lambda \boxtimes P)$ to $I \cup \{J\}$ (where $J$ represents no matching), followed by type changing for the matched agents that is modeled by a process $\alpha^n$ from $(I \times \Omega, \mathcal{I} \boxtimes \mathcal{F}, \lambda \boxtimes P)$ to $S$.

For the random mutation step at time $n$, given a $K \times K$ probability transition matrix[13] $b$, we require that, for each agent $i \in I$,

$$P(h_i^n = l | \alpha_i^{n-1} = k) = b_{kl}, \tag{1}$$

the specified probability with which an agent $i$ of type $k$ at the end of time period $n-1$ mutates to type $l$.

For the random partial matching step at time $n$, we let $\overline{p}^{n-1/2}$ be the expected cross-sectional type distribution immediately after random mutation. That is,

$$\overline{p}_k^{n-1/2} = \overline{p}^{n-1/2}(k) = \int_\Omega \lambda(\{i \in I : h_\omega^n(i) = k\}) \, dP(\omega). \tag{2}$$

The random partial matching function $\pi^n$ at time $n$ is defined by:

1. For any $\omega \in \Omega$, $\pi_\omega^n(\cdot)$ is a full matching on $I - (\pi_\omega^n)^{-1}(\{J\})$, as defined in Section 2.3.
2. Extending $h^n$ so that $h^n(J, \omega) = J$ for any $\omega \in \Omega$, we define $g^n : I \times \Omega \to S \cup \{J\}$ by

$$g^n(i, \omega) = h^n(\pi^n(i, \omega), \omega),$$

and assume that $g^n$ is $\mathcal{I} \boxtimes \mathcal{F}$-measurable.
3. Let $q \in [0,1]^S$. For each agent $i \in I$,

$$P(g_i^n = J | h_i^n = k) = q_k,$$

$$P(g_i^n = l | h_i^n = k) = \frac{(1-q_k)(1-q_l)\overline{p}_l^{n-1/2}}{\sum_{r=1}^K (1-q_r)\overline{p}_r^{n-1/2}}. \tag{3}$$

Equation (3) means that for any agent whose type before the matching is $k$, the probability of being unmatched is $q_k$, and the probability of being matched to a type-$l$ agent is proportional to the expected cross-sectional type distribution for matched agents. When $g^n$ is essentially pairwise independent (as under the Markov conditional independence condition in Section 3.2 below), the exact law of large numbers in [38] and [40] implies that the realized cross-sectional type distribution $\lambda(h_\omega^n)^{-1}$ after random mutation at time $n$ is indeed the expected distribution $\overline{p}^{n-1/2}$, $P$-almost surely.[14]

---

[13]Here, $b_{kl}$ is in $[0,1]$, with $\sum_{l=1}^K b_{kl} = 1$ for each $k$.

[14]As noted in Footnote 11, if the denominator in (3) is zero, then almost no agents will be matched and we can simply interpret the ratio as zero.

Finally, for the step of random type changing for matched agents at time $n$, a given $\nu\colon S\times S\to\Delta$ specifies the probability distribution $\nu_{kl}=\nu(k,l)$ of the new type of a type-$k$ agent who has met a type-$l$ agent. When agent $i$ is not matched at time $n$, she keeps her type $h_i^n$ with probability 1. We thus require that the type function $\alpha^n$ after matching satisfies, for each agent $i\in I$,

$$P(\alpha_i^n = r | h_i^n = k, g_i^n = J) = \delta_k^r,$$
(4)
$$P(\alpha_i^n = r | h_i^n = k, g_i^n = l) = \nu_{kl}(r),$$

where $\delta_k^r$ is 1 if $r=k$, and zero otherwise.

Thus, we have inductively defined a dynamical system $\mathbb{D}$ with random mutation, partial matching and match-induced type changing with parameters $(p^0, b, q, \nu)$.

3.2. *Markov conditional independence in types.* In this subsection, we consider a suitable independence condition on the dynamical system $\mathbb{D}$. For $n\geq 1$, to formalize the intuitive idea that given their type function $\alpha^{n-1}$, the agents randomly mutate to other types independently at time $n$, and that their types in earlier periods have no effect on this mutation, we say that the random mutation is *Markov conditionally independent in types* if, for $\lambda$-almost all $i\in I$ and $\lambda$-almost all $j\in I$,

(5)
$$P(h_i^n = k, h_j^n = l | \alpha_i^0, \ldots, \alpha_i^{n-1}; \alpha_j^0, \ldots, \alpha_j^{n-1})$$
$$= P(h_i^n = k | \alpha_i^{n-1}) P(h_j^n = l | \alpha_j^{n-1})$$

holds for all types $k,l\in S$.[15]

Intuitively, the random partial matching at time $n$ should depend only on agents' types immediately after the random mutation. One may also want the random partial matching to be independent across agents, given events that occurred in the first $n-1$ time periods and the random mutation at time $n$. We say that the random partial matching $\pi^n$ is *Markov conditionally independent in types* if, for $\lambda$-almost all $i\in I$ and $\lambda$-almost all $j\in I$,

(6)
$$P(g_i^n = c, g_j^n = d | \alpha_i^0, \ldots, \alpha_i^{n-1}, h_i^n; \alpha_j^0, \ldots, \alpha_j^{n-1}, h_j^n)$$
$$= P(g_i^n = c | h_i^n) P(g_j^n = d | h_j^n)$$

holds for all types $c,d\in S\cup\{J\}$.

The agents' types at the end of time period $n$ should depend on the agents' types immediately after the random mutation stage at time $n$, as well as the results of random partial matching at time $n$, but not otherwise on events that occurred in previous periods. This motivates the following definition. The random type changing after partial matching at time $n$ is said to be *Markov conditionally independent in types* if for $\lambda$-almost all $i \in I$ and $\lambda$-almost all $j \in I$, and for each $n \geq 1$,

$$
\begin{aligned}
(7) \quad & P(\alpha_i^n = k, \alpha_j^n = l | \alpha_i^0, \ldots, \alpha_i^{n-1}, h_i^n, g_i^n; \alpha_j^0, \ldots, \alpha_j^{n-1}, h_j^n, g_j^n) \\
& = P(\alpha_i^n = k | h_i^n, g_i^n) P(\alpha_j^n = l | h_j^n, g_j^n)
\end{aligned}
$$

holds for all types $k, l \in S$.

---

[15] We could include the functions $h^m$ and $g^m$ for $1\leq m\leq n-1$ as well. However, it is not necessary to do so since we only care about the dependence structure across time for the type functions at the end of each time period.

The dynamical system $\mathbb{D}$ is said to be *Markov conditionally independent in types* if, in each time period $n$, each random step (random mutation, partial matching and type changing) is so.

3.3. *The existence theorem.* The following theorem shows the existence of a dynamical system with random mutation, partial matching and type changing that is Markov conditionally independent in types. Its proof will be given in Section 6 after the development of a generalized Fubini theorem for a Loeb transition probability and a Loeb product transition probability system in Section 5.

THEOREM 3.1. *Fixing any parameters $p^0$ for the initial cross-sectional type distribution, $b$ for mutation probabilities, $q \in [0,1]^S$ for no-match probabilities and $\nu$ for match-induced type-change probabilities, there exists a Fubini extension of the usual product probability space on which is defined a dynamical system $\mathbb{D}$ with random mutation, partial matching and type changing that is Markov conditionally independent in types with these parameters $(p^0, b, q, \nu)$.*

Note that the dynamic matching model $\mathbb{D}$ described in Theorem 3.1 above satisfies the conditions in Theorem 3 of [15]. Thus, the conclusions of Theorem 3 in [15] also hold for the matching model $\mathbb{D}$, including the statements that the type processes of individual agents in such a dynamical system $\mathbb{D}$ form a continuum of independent Markov chains, and that the time evolution of the cross-sectional type process is completely determined from a Markov chain with known transition matrices.

**4. Proofs of Theorems 2.4 and 2.6.** In this section, we give the proofs for Theorems 2.4 and 2.6 in Sections 4.1 and 4.2, respectively.

The space of agents used from this section onward will be based on a hyperfinite Loeb counting probability space $(I, \mathcal{I}, \lambda)$ that is the Loeb space (see [33] and [34]) of internal probability space $(I, \mathcal{I}_0, \lambda_0)$, where $I$ is a hyperfinite set, $\mathcal{I}_0$ its internal power set, and $\lambda_0(A) = |A|/|I|$ for any $A \in \mathcal{I}_0$ (i.e., $\lambda_0$ is the internal counting probability measure on $I$). Using the usual ultrapower construction as in [34], the hyperfinite set $I$ itself can be viewed as an equivalence class of a sequence of finite sets whose sizes go to infinity, and the external cardinality of $I$ is the cardinality of the continuum. Thus, one can also take the unit interval $[0,1]$ as the space of agents, endowed with a $\sigma$-algebra and an atomless probability measure via a bijection between $I$ and the unit interval.

All of the internal probability spaces to be discussed from this section onward are hyperfinite internal probability spaces. A general hyperfinite internal probability space is an ordered triple $(\Omega, \mathcal{F}_0, P_0)$, where $\Omega = \{\omega_1, \omega_2, \ldots, \omega_\gamma\}$ for some unlimited hyperfinite natural number $\gamma$, $\mathcal{F}_0$ is the internal power set on $\Omega$, and $P_0(B) = \sum_{1 \leq j \leq \gamma, \omega_j \in B} P_0(\{\omega_j\})$ for any $B \in \mathcal{F}_0$. When the weights $P_0(\{\omega_j\})$, $1 \leq j \leq \gamma$, are all infinitesimals, $(\Omega, \mathcal{F}_0, P_0)$ is said to be atomless, and its Loeb space $(\Omega, \mathcal{F}, P)$, as a standard probability space, is atomless in the usual sense of Footnote 8. Note that nonstandard analysis is used extensively from this section onward. The reader is referred to the first three chapters of [34] for more details.

4.1. *Proof of Theorem* 2.4. Fix an even hyperfinite natural number in the set $^*\mathbb{N}_\infty$ of unlimited hyperfinite natural numbers. Let $I = \{1, 2, \ldots, N\}$, let $\mathcal{I}_0$ be the collection of all the internal subsets of $I$, and let $\lambda_0$ be the internal counting probability measure on $\mathcal{I}_0$. Let $(I, \mathcal{I}, \lambda)$ be the Loeb space of the internal probability space $(I, \mathcal{I}_0, \lambda_0)$. Note that $(I, \mathcal{I}, \lambda)$ is obviously atomless.

We can draw agents from $I$ in pairs without replacement, and then match them in these pairs. The procedure can be the following. Take one fixed agent; this agent can be matched with $N-1$ different agents. After the first pair is matched, there are $N-2$ agents. We can do the same thing to match a second pair with $N-3$ possibilities. Continue this procedure to produce a total number of $1 \times 3 \times \cdots \times (N-3) \times (N-1)$, denoted by $(N-1)!!$, different matchings. Let $\Omega$ be the space of all such matchings, $\mathcal{F}_0$ the collection of all internal subsets of $\Omega$ and $P_0$ the internal counting probability measure on $\mathcal{F}_0$. Let $(\Omega, \mathcal{F}, P)$ be the Loeb space of the internal probability space $(\Omega, \mathcal{F}_0, P_0)$.

Let $(I \times \Omega, \mathcal{I}_0 \otimes \mathcal{F}_0, \lambda_0 \otimes P_0)$ be the internal product probability space of $(I, \mathcal{I}_0, \lambda_0)$ and $(\Omega, \mathcal{F}_0, P_0)$. Then $\mathcal{I}_0 \otimes \mathcal{F}_0$ is actually the collection of all the internal subsets of $I \times \Omega$ and $\lambda_0 \otimes P_0$ is the internal counting probability measure on $\mathcal{I}_0 \otimes \mathcal{F}_0$. Let $(I \times \Omega, \mathcal{I} \boxtimes \mathcal{F}, \lambda \boxtimes P)$ be the Loeb space of the internal product $(I \times \Omega, \mathcal{I}_0 \otimes \mathcal{F}_0, \lambda_0 \otimes P_0)$, which is indeed a Fubini extension of the usual product probability space.[16]

---

[16] For any given two Loeb spaces $(I, \mathcal{I}, \lambda)$ and $(\Omega, \mathcal{F}, P)$, it is shown in [30] that the Loeb product space $(I \times \Omega, \mathcal{I} \boxtimes \mathcal{F}, \lambda \boxtimes P)$ is uniquely defined by the marginal Loeb spaces.

Now, for a given matching $\omega \in \Omega$ and a given agent $i$, let $\pi(i,\omega)$ be the unique $j$ such that the pair $(i,j)$ is matched under $\omega$. For each $\omega \in \Omega$, since $\pi_\omega$ is an internal bijection on $I$, it is obvious that $\pi_\omega$ is measure-preserving from the Loeb space $(I, \mathcal{I}, \lambda)$ to itself. Thus, (i) of part (1) is shown.

It is obvious that for any agent $i \in I$,

$$(8) \qquad P_0(\{\omega \in \Omega : \pi_i(\omega) = j\}) = \frac{1}{N-1}$$

for any $j \neq i$; that is, the $i$th agent is matched with equal chance to other agents.

Fix any $i \in I$. For any internal set $C \in \mathcal{I}_0$, (8) implies that $P_0(\omega \in \Omega : \pi_i(\omega) \in C)$ is $|C|/(N-1)$ if $i \notin C$, and $(|C|-1)/(N-1)$ if $i \in C$, where $|C|$ is the internal cardinality of $C$. This means[17] that

$$(9) \qquad P_0(\pi_i^{-1}(C)) \simeq \frac{|C|}{N} = \lambda_0(C) \simeq \lambda(C).$$

Therefore, $\pi_i$ is a measure-preserving mapping from $(\Omega, \mathcal{F}_0, P)$ to $(I, \mathcal{I}_0, \lambda)$, and is measure-preserving from the Loeb space $(\Omega, \mathcal{F}, P)$ to the Loeb space $(I, \mathcal{I}, \lambda)$.[18] Thus, (ii) of part (1) is shown.

We can also obtain that

$$(10) \qquad (\lambda_0 \otimes P_0)(\pi^{-1}(C)) = \int_{i \in I} P_0(\pi_i^{-1}(C)) \, d\lambda_0(i) \simeq \lambda_0(C) \simeq \lambda(C).$$

A proof similar to that of Footnote 18 shows that $\pi$ is a measure-preserving mapping from $(I \times \Omega, \mathcal{I} \boxtimes \mathcal{F}, \lambda \boxtimes P)$ to $(I, \mathcal{I}, \lambda)$.

Next, for $i \neq j$, consider the joint event

$$(11) \qquad E = \{\omega \in \Omega : (\pi_i(\omega), \pi_j(\omega)) = (i', j')\},$$

---

Anderson noted in [3] that $(I \times \Omega, \mathcal{I} \boxtimes \mathcal{F}, \lambda \boxtimes P)$ is an extension of the usual product $(I \times \Omega, \mathcal{I} \otimes \mathcal{F}, \lambda \otimes P)$. Keisler proved in [26] (see also [34]) that the Fubini property still holds on $(I \times \Omega, \mathcal{I} \boxtimes \mathcal{F}, \lambda \boxtimes P)$. Thus, the Loeb product space is a Fubini extension of the usual product probability space. In addition, it is shown in Theorem 6.2 of [38] that when both $\lambda$ and $P$ are atomless, $(I \times \Omega, \mathcal{I} \boxtimes \mathcal{F}, \lambda \boxtimes P)$ is rich enough to be endowed with a process $h$ whose random variables are essentially pairwise independent and can take any variety of distributions (and in particular the uniform distribution on $[0,1]$).

[17]For two hyperreals $\alpha$ and $\beta$, $\alpha \simeq \beta$ means that the difference $\alpha - \beta$ is an infinitesimal; see [34].

[18]For any Loeb measurable set $B \in \mathcal{I}$ and for any standard positive real number $\varepsilon$, there are internal sets $C$ and $D$ in $\mathcal{I}_0$ such that $C \subseteq B \subseteq D$ and $\lambda_0(D - C) < \varepsilon$. Thus $\pi_i^{-1}(C) \subseteq \pi_i^{-1}(B) \subseteq \pi_i^{-1}(D)$, and

$$P_0(\pi_i^{-1}(D) - \pi_i^{-1}(C)) \simeq \lambda_0(D - C) < \varepsilon,$$

which implies that $\pi_i^{-1}(B)$ is Loeb measurable in $\mathcal{F}$. Also, $\lambda(C) \leq P(\pi_i^{-1}(B)) \leq \lambda(D)$, and thus $|P(\pi_i^{-1}(B)) - \lambda(B)| \leq \lambda(D - C) \leq \varepsilon$ for any standard positive real number $\varepsilon$. This means that $P(\pi_i^{-1}(B)) = \lambda(B)$.

that is, the $i$th agent is matched to the $i'$th agent and the $j$th agent is matched to the $j'$th agent. In order to show the measure-preserving property of the mapping $(\pi_i, \pi_j)$ in the following paragraph, we need to know the value of $P_0(E)$ in three different cases. The first case is $(i', j') = (j, i)$, that is, the $i$th agent is matched to the $j$th agent and the $j$th agent is matched to the $i$th agent. In this case, $P_0(E) = 1/(N-1)$. The second case is that of $P_0(E) = 0$, which holds when $i' = i$ or $j' = j$ (the $i$th agent is matched to herself, or the $j$th agent is matched to herself), or when $i' = j$ but $j' \neq i$ (the $i$th agent is matched to the $j$th agent, but the $j$th agent is *not* matched to the $i$th agent), or when $j' = i$ but $i' \neq j$ (the $j$th agent is matched to the $i$th agent, but the $i$th agent is *not* matched to the $j$th agent), or when $i' = j'$ (both the $i$th agent and the $j$th agent are matched to the same agent). The third case applies if the indices $i, j$ and $i', j'$ are completely distinct. In this third case, after the pairs $(i, i')$, $(j, j')$ are drawn, there are $N-4$ agents left, and hence there are $(N-5)!!$ ways to draw the rest of the pairs in order to complete the matching. This means that $P_0(E) = (N-5)!!/(N-1)!! = 1/((N-1)(N-3))$.

Let $(I \times I, \mathcal{I}_0 \otimes \mathcal{I}_0, \lambda_0 \otimes \lambda_0)$ be the internal product of $(I, \mathcal{I}_0, \lambda_0)$ with itself, and $(I \times I, \mathcal{I} \boxtimes \mathcal{I}, \lambda \boxtimes \lambda)$ the Loeb space of the internal product. Fix any $i, j \in I$ with $i \neq j$. Let $D$ be the diagonal $\{(i', i') : i' \in I\}$. The third case of the above paragraph implies that for any internal set $G \in \mathcal{I}_0 \otimes \mathcal{I}_0$,

$$\begin{aligned}
(12) \quad & P_0(\{\omega \in \Omega : (\pi_i(\omega), \pi_j(\omega)) \in G - (D \cup (\{i,j\} \times I) \cup (I \times \{i,j\}))\}) \\
& = \frac{|G - (D \cup (\{i,j\} \times I) \cup (I \times \{i,j\}))|}{(N-1)(N-3)} \simeq \frac{|G|}{(N)^2} = (\lambda_0 \otimes \lambda_0)(G).
\end{aligned}$$

By using the formula for $P_0(E)$ in the first two cases, we can obtain that

$$\begin{aligned}
(13) \quad & P_0(\{\omega \in \Omega : (\pi_i(\omega), \pi_j(\omega)) \in (D \cup (\{i,j\} \times I) \cup (I \times \{i,j\}))\}) \\
& = \frac{1}{N-1} \simeq 0.
\end{aligned}$$

Equations (12) and (13) imply that

$$(14) \qquad P_0(\{\omega \in \Omega : (\pi_i(\omega), \pi_j(\omega)) \in G\}) \simeq (\lambda_0 \otimes \lambda_0)(G).$$

A proof similar to that of Footnote 18 shows that $(\pi_i, \pi_j)$ is a measure-preserving mapping from $(\Omega, \mathcal{F}, P)$ to $(I \times I, \mathcal{I} \boxtimes \mathcal{I}, \lambda \boxtimes \lambda)$.

Let $\alpha$ be an $\mathcal{I}$-measurable type function with a distribution $p$ on $S$, $I_k = \alpha^{-1}(\{k\})$ and $p_k = \lambda(I_k)$ for $1 \leq k \leq K$. Let $g = \alpha(\pi)$. Then, for any $1 \leq k \leq K$, $g^{-1}(\{k\}) = \pi^{-1}(I_k)$, which is Loeb product measurable in $\mathcal{I} \boxtimes \mathcal{F}$ with $\lambda \boxtimes P$-measure $p_k$ because of the measure-preserving property of $\pi$. Hence, $g$ is $\mathcal{I} \boxtimes \mathcal{F}$-measurable. For each $i \in I$, the measure-preserving property of $\pi_i$ implies that $g_i$ has the same distribution $p$ as $\alpha$.

Fix any $i, j \in I$ with $i \neq j$. For any $1 \leq k, l \leq K$, the measure-preserving property of $(\pi_i, \pi_j)$ implies that

$$
\begin{aligned}
P(g_i = k, g_j = l) &= P(\{\omega \in \Omega : (\pi_i(\omega), \pi_j(\omega)) \in I_k \times I_l\}) \\
&= (\lambda \boxtimes \lambda)(I_k \times I_l) = P(g_i = k) \cdot P(g_j = l),
\end{aligned}
\tag{15}
$$

which means that the random variables $g_i$ and $g_j$ are independent.[19] Hence, part (2) is shown.

Finally, take any $A_1, A_2 \in \mathcal{I}$, and let $f(i, \omega) = 1_{A_2}(\pi(i, \omega))$ for all $(i, \omega) \in I \times \Omega$. Then, $f$ is an i.i.d. process with a common distribution on $\{0, 1\}$ with probabilities $\lambda(A_2)$ on $\{1\}$ and $1 - \lambda(A_2)$ on $\{0\}$. By the exact law of large numbers in Theorem 5.2 of [38] or Theorem 3.5 of [40],[20] one has for $P$-almost all $\omega \in \Omega$, $\lambda^{A_1}(f^{A_1})_\omega^{-1}(\{1\}) = \lambda(A_2)$, which implies that $\lambda(A_2) = \lambda(A_1 \cap \pi_\omega^{-1}(A_2))/\lambda(A_1)$. Hence, (iii) of part (1) follows.

4.2. *Proof of Theorem* 2.6. Let $M$ be any fixed unlimited hyperfinite natural number in $^*\mathbb{N}_\infty$, and let $I = \{1, 2, \ldots, M\}$ be the space of agents. Let $\mathcal{I}_0$ be the collection of all the internal subsets of $I$, and $\lambda_0$ the internal counting probability measure on $\mathcal{I}_0$. Let $(I, \mathcal{I}, \lambda)$ be the Loeb space of the internal probability space $(I, \mathcal{I}_0, \lambda_0)$.

Let $\beta$ be any $\mathcal{I}$-measurable type function from $I$ to $S = \{1, \ldots, K\}$. We can find an internal type function $\alpha$ from $I$ to $S$ such that $\lambda(\{i \in I : \alpha(i) \neq \beta(i)\}) = 0$. Let $A_k = \alpha^{-1}(k)$ and $|A_k| = M_k$ for $1 \leq k \leq K$ with $\sum_{k=1}^K M_k = M$. Then $\lambda(A_k) = p_k \simeq M_k/M$ for $1 \leq k \leq K$. Without loss of generality, we can assume that $M_k \in {}^*\mathbb{N}_\infty$.[21] For each $k$ in $\{1, \ldots, K\}$, pick a hyperfinite natural number $m_k$ such that $M_k - m_k \in {}^*\mathbb{N}_\infty$ and $q_k \simeq m_k/M_k$, and such that

---

[19]When the type function $\alpha$ is internal, the type process is internal as well. However, the computations in (12)–(14) indicate that the random variables $g_i$, $i \in I$, are in general not $*$-independent; see, for example, [3] for the definition of $*$-independence.

[20]What we need in this paper is a special case of Theorem 3.5 in [40]. Let $f$ be a process from $(I \times \Omega, \mathcal{I} \boxtimes \mathcal{F}, \lambda \boxtimes P)$ to a complete separable metric space $X$. Assume that the random variables $f_i$, $i \in I$ are pairwise independent. Then, for $P$-almost all $\omega \in \Omega$, the sample function $f_\omega$ has the same distribution as $f$ in the sense that for any Borel set $B$ in $X$, $\lambda(f_\omega^{-1}(B)) = (\lambda \boxtimes P)(f^{-1}(B))$. Fix any $A \in \mathcal{I}$ with $\lambda(A) > 0$. Let $f^A$ be the restriction of $f$ to $A \times \Omega$, $\lambda^A$ and $\lambda^A \boxtimes P$ the probability measures rescaled from the restrictions of $\lambda$ and $\lambda \boxtimes P$ to $\{D \in \mathcal{I} : D \subseteq A\}$ and $\{C \in \mathcal{I} \boxtimes \mathcal{F} : C \subseteq A \times \Omega\}$, respectively. Then, for the case that the random variables $f_i, i \in I$ have a common distribution $\mu$ on $X$, the sample function $(f^A)_\omega$ also has distribution $\mu$ for $P$-almost all $\omega \in \Omega$.

[21]When $p_k = 0$, we may still need to divide some number $m_k$ by $M_k$ so that the ratio is infinitely close to a real number $q_k$. We can take $M_k \in {}^*\mathbb{N}_\infty$ with $M_k/M \simeq 0$. We then take an internal subset of $I$ with $M_k$ many points as $A_k$ and adjust the rest $A_l, l \neq k$, on some $\lambda$-null internal sets. This will produce a new internal type function $\alpha$ with the desired properties.

$N = \sum_{l=1}^{K}(M_l - m_l)$ is an unlimited even hyperfinite natural number. It is easy to see that

(16) $$\frac{N}{M} = \sum_{l=1}^{K} \frac{M_l}{M}\left(1 - \frac{m_l}{M_l}\right) \simeq \sum_{l=1}^{K} p_l(1 - q_l).$$

For each $k$ in $\{1, 2, \ldots, K\}$, let $B_k$ be an arbitrary internal subset of $A_k$ with $m_k$ elements, and let $\mathcal{P}_{m_k}(A_k)$ be the collection of all such internal subsets. For given $B_k \in \mathcal{P}_{m_k}(A_k)$ for $k = 1, 2, \ldots, K$, let $\pi^{B_1, B_2, \ldots, B_K}$ be a (full) matching on $I - \bigcup_{k=1}^{K} B_k$ produced by the method described in the proof of Theorem 2.4; there are $(N-1)!! = 1 \times 3 \times \cdots \times (N-3) \times (N-1)$ such matchings.

Our sample space $\Omega$ is the set of all ordered tuples $(B_1, B_2, \ldots, B_K, \pi^{B_1, B_2, \ldots, B_K})$ such that $B_k \in \mathcal{P}_{m_k}(A_k)$ for each $k = 1, \ldots, K$, and $\pi^{B_1, B_2, \ldots, B_K}$ is a matching on $I - \bigcup_{k=1}^{K} B_k$. Then, $\Omega$ has $((N-1)!!) \prod_{k=1}^{K} \binom{M_k}{m_k}$ many elements in total. Let $P_0$ be the internal counting probability measure defined on the collection $\mathcal{F}_0$ of all the internal subsets of $\Omega$. Let $(\Omega, \mathcal{F}, P)$ be the Loeb space of the internal probability space $(\Omega, \mathcal{F}_0, P_0)$. Note that both $(I, \mathcal{I}_0, \lambda_0)$ and $(\Omega, \mathcal{F}_0, P_0)$ are atomless. Let $(I \times \Omega, \mathcal{I} \boxtimes \mathcal{F}, \lambda \boxtimes P)$ be the Loeb space of the internal product $(I \times \Omega, \mathcal{I}_0 \otimes \mathcal{F}_0, \lambda_0 \otimes P_0)$, which is a Fubini extension of the usual product probability space by Footnote 16.

Let $J$ represent nonmatching. Define a mapping $\pi$ from $I \times \Omega$ to $I \cup \{J\}$. For $i \in A_k$ and $\omega = (B_1, B_2, \ldots, B_K, \pi^{B_1, B_2, \ldots, B_K})$, if $i \in B_k$, then $\pi(i, \omega) = J$ (agent $i$ is not matched); if $i \notin B_k$, then $i \in I - \bigcup_{r=1}^{K} B_r$, agent $i$ is to be matched with agent $\pi^{B_1, B_2, \ldots, B_K}(i)$, and let $\pi(i, \omega) = \pi^{B_1, B_2, \ldots, B_K}(i)$. It is obvious that $\pi_\omega^{-1}(\{J\}) = \bigcup_{r=1}^{K} B_r$ and that the restriction of $\pi_\omega$ to $I - \bigcup_{r=1}^{K} B_r$ is a full matching on the set. Let $g^\beta$ be the matched type process from $I \times \Omega$ to $S \cup \{J\}$ under the type function $\beta$; that is, $g^\beta(i, \omega) = \beta(\pi(i, \omega))$ with $\beta(J) = J$.

When $\sum_{r=1}^{K} p_r(1 - q_r) = 0$, we know that $N/M \simeq 0$. For those $i \in A_k$ with $p_k > 0$, it is clear that $P_0(\{\omega \in \Omega : \pi_i(\omega) = J\}) = m_k/M_k \simeq q_k = 1$, and thus $\lambda \boxtimes P(\pi(i, \omega) \neq J) = 0$, which means that $\lambda \boxtimes P(g^\beta(i, \omega) \neq J) = 0$, and $g_i^\beta(\omega) = J$ for $P$-almost all $\omega \in \Omega$. Thus conditions (1) and (2) in Definition 2.5 are satisfied trivially; that is, one has a trivial random partial matching that is independent in types.

For the rest of the proof, assume that $\sum_{r=1}^{K} p_r(1 - q_r) > 0$. Let $g$ be the matched type process from $I \times \Omega$ to $S \cup \{J\}$, defined by $g(i, \omega) = \alpha(\pi(i, \omega))$, where $\alpha(J) = J$. Since both $\alpha$ and $\pi$ are internal, the fact that $\mathcal{I}_0 \otimes \mathcal{F}_0$ is the internal power set on $I \times \Omega$ implies that $g$ is $\mathcal{I}_0 \otimes \mathcal{F}_0$-measurable, and thus $\mathcal{I} \boxtimes \mathcal{F}$-measurable.

Fix an agent $i \in A_k$ for some $1 \leq k \leq K$. For any $1 \leq l \leq K$, and for any $B_r \in \mathcal{P}_{m_r}(A_r)$, $r = 1, 2, \ldots, K$, let $N_{il}^{B_1, B_2, \ldots, B_K}$ be the number of full matchings on $\bigcup_{r=1}^{K}(A_r - B_r)$ such that agent $i$ is matched to some agent in $A_l - B_l$. It is obvious that $N_{il}^{B_1,B_2,\ldots,B_K}$ depends only on the numbers of points in the sets $A_r - B_r$, $r = 1,\ldots,K$, which are $M_r - m_r$, $r = 1,\ldots,K$, respectively. Hence, $N_{il}^{B_1,B_2,\ldots,B_K}$ is independent of the particular choices of $B_1, B_2, \ldots, B_K$, and so can simply be denoted by $N_{il}$. Then, (9) implies that

$$\frac{N_{il}}{(N-1)!!} \simeq \frac{M_l - m_l}{N}. \tag{17}$$

It can be checked that the internal cardinality of the event $\{g_i = l\}$ is

$$|\{\omega \in \Omega : i \in A_k - B_k, \pi_i^{B_1,B_2,\ldots,B_K}(\omega) \in (A_l - B_l)\}|$$
$$= \binom{M_k - 1}{m_k} \left(\prod_{r \neq k} \binom{M_r}{m_r}\right) N_{il}, \tag{18}$$

for $\omega = (B_1, B_2, \ldots, B_K, \pi^{B_1,B_2,\ldots,B_K})$. Hence (17) and (18) imply that

$$P_0(g_i = l) = \frac{M_k - m_k}{M_k} \frac{N_{il}}{(N-1)!!} \simeq (1 - q_k)\frac{M_l - m_l}{N}$$
$$= (1 - q_k)\frac{(1 - m_l/M_l)M_l/M}{(N/M)} \simeq \frac{(1 - q_k)p_l(1 - q_l)}{\sum_{r=1}^K p_r(1 - q_r)}. \tag{19}$$

It is also easy to see that $P_0(\{\omega \in \Omega : g_i(\omega) = J\}) = m_k/M_k \simeq q_k$. This means that for $i \in A_k$,

$$P(g_i = l) = \frac{(1 - q_k)p_l(1 - q_l)}{\sum_{r=1}^K p_r(1 - q_r)},$$

and that $P(g_i = J) = q_k$. Hence, the distribution condition on $g_i$ is satisfied for each $i \in I$.

We need to show that the random partial matching $\pi$ is independent in types. Fix agents $i, j \in I$ with $i \neq j$. For any $1 \leq l, t \leq K$, and for any $B_r \in \mathcal{P}_{m_r}(A_r)$, $r = 1, 2, \ldots, K$, let $N_{iljt}^{B_1,B_2,\ldots,B_K}$ be the number of full matchings on $\bigcup_{r=1}^K (A_r - B_r)$ such that agents $i$ and $j$ are matched to some agents respectively in $A_l - B_l$ and $A_t - B_t$. As in the case of $N_{il}^{B_1,B_2,\ldots,B_K}$, $N_{iljt}^{B_1,B_2,\ldots,B_K}$ is independent of the particular choices of $B_1, B_2, \ldots, B_K$ and can simply be denoted by $N_{iljt}$. By taking $G = (A_l - B_l) \times (A_t - B_t)$, (14) implies that

$$\frac{N_{iljt}}{(N-1)!!} \simeq \frac{M_l - m_l}{N} \frac{M_t - m_t}{N}. \tag{20}$$

We first consider the case that both $i$ and $j$ belong to $A_k$ for some $k$ in $\{1, \ldots, K\}$. It is easy to see that $P_0(g_i = J, g_j = J) = m_k(m_k - 1)/(M_k(M_k - 1))$, and hence that

$$P(g_i = J, g_j = J) = P(g_i = J)P(g_j = J) = q_k^2. \tag{21}$$

As above, it can be checked that the internal cardinality of the event $\{g_i = l, g_j = J\}$ is

$$
(22) \quad |\{\omega \in \Omega : i \in A_k - B_k, j \in B_k, \pi_i^{B_1,\ldots,B_K}(\omega) \in (A_l - B_l)\}|
$$
$$
= \binom{M_k - 2}{m_k - 1} \left(\prod_{r \neq k} \binom{M_r}{m_r}\right) N_{il},
$$

for $\omega = (B_1, \ldots, B_K, \pi^{B_1,\ldots,B_K})$. Hence (17) and (22) imply that

$$
P_0(g_i = l, g_j = J) = \frac{m_k(M_k - m_k)}{M_k(M_k - 1)} \frac{N_{il}}{(N-1)!!}
$$
$$
(23) \quad \simeq q_k(1 - q_k)\frac{M_l - m_l}{N}
$$
$$
\simeq \frac{q_k(1 - q_k)p_l(1 - q_l)}{\sum_{r=1}^{K} p_r(1 - q_r)},
$$

which implies that

$$
(24) \quad P(g_i = l, g_j = J) = P(g_i = l)P(g_j = J).
$$

Similarly, the events $(g_i = J)$ and $(g_j = l)$ are independent.

The event $\{g_i = l, g_j = t\}$ is actually the set of all the $\omega = (B_1, \ldots, B_K, \pi^{B_1,\ldots,B_K})$ such that both $i$ and $j$ are in $A_k - B_k$, and agents $i$ and $j$ are matched to some agents in $A_l - B_l$ and $A_t - B_t$, respectively. Thus, the internal cardinality of $\{g_i = l, g_j = t\}$ is

$$
(25) \quad \binom{M_k - 2}{m_k} \left(\prod_{r \neq k} \binom{M_r}{m_r}\right) N_{iljt}.
$$

Hence (20) and (25) imply that

$$
P_0(g_i = l, g_j = t) = \frac{(M_k - m_k)(M_k - m_k - 1)}{M_k(M_k - 1)} \frac{N_{iljt}}{(N-1)!!}
$$
$$
(26) \quad \simeq (1 - q_k)^2 \frac{M_l - m_l}{N} \frac{M_t - m_t}{N}
$$
$$
\simeq \frac{(1 - q_k)^2 p_l(1 - q_l)p_t(1 - q_t)}{(\sum_{r=1}^{K} p_r(1 - q_r))^2},
$$

which implies that

$$
(27) \quad P(g_i = l, g_j = t) = P(g_i = l)P(g_j = t).
$$

Hence the random variables $g_i$ and $g_j$ are independent.

For the case that $i \in A_k$ and $j \in A_n$ with $1 \leq k \neq n \leq K$, one can first observe that

(28)
$$P_0(g_i = l, g_j = J) = \frac{(M_k - m_k)m_n}{M_k M_n} \frac{N_{il}}{(N-1)!!},$$
$$P_0(g_i = l, g_j = t) = \frac{(M_k - m_k)(M_n - m_n)}{M_k M_n} \frac{N_{iljt}}{(N-1)!!}.$$

In this case, one can use computations similar to those of the above two paragraphs to show that the random variables $g_i$ and $g_j$ are independent. The details are omitted here.

We have proved the result for the type function $\alpha$. We still need to prove it for $\beta$ [and for $g^\beta = \beta(\pi)$]. Fix any agent $i \in A_k$, for some $1 \leq k \leq K$. For any internal set $A \in \mathcal{I}_0$, and for any $B_r \in \mathcal{P}_{m_r}(A_r)$, $r = 1, 2, \ldots, K$, let $N_{iA}^{B_1, B_2, \ldots, B_K}$ be the number of full matchings on $\bigcup_{r=1}^{K}(A_r - B_r)$ such that agent $i$ is matched to some agent in $A - \bigcup_{r=1}^{K} B_r$. Then, (9) implies that

(29)
$$\frac{N_{iA}^{B_1, B_2, \ldots, B_K}}{(N-1)!!} \simeq \frac{|A - \bigcup_{r=1}^{K} B_r|}{N}.$$

The internal event $\pi_i^{-1}(A)$ is, for $\omega = (B_1, B_2, \ldots, B_K, \pi^{B_1, B_2, \ldots, B_K})$,

(30)
$$\left\{ \omega \in \Omega : i \in A_k - B_k, \pi_i^{B_1, B_2, \ldots, B_K}(\omega) \in \left( A - \bigcup_{r=1}^{K} B_r \right) \right\}.$$

Hence (29) and (30) imply that

(31)
$$P_0(\pi_i^{-1}(A)) = \sum_{B_k \in \mathcal{P}_{m_k}(A_k \setminus \{i\}), B_l \in \mathcal{P}_{m_l}(A_l) \text{ for } l \neq k} \frac{1}{\prod_{r=1}^{K} \binom{M_r}{m_r}} \frac{N_{iA}^{B_1, B_2, \ldots, B_K}}{(N-1)!!}$$
$$\simeq \sum_{B_k \in \mathcal{P}_{m_k}(A_k \setminus \{i\}), B_l \in \mathcal{P}_{m_l}(A_l) \text{ for } l \neq k} \frac{1}{\prod_{r=1}^{K} \binom{M_r}{m_r}} \frac{|A - \bigcup_{r=1}^{K} B_r|}{N}$$
$$\leq \frac{|A|}{N} \frac{\binom{M_k - 1}{m_k}(\prod_{r \neq k} \binom{M_r}{m_r})}{\prod_{r=1}^{K} \binom{M_r}{m_r}} = \frac{M_k - m_k}{M_k} \frac{|A|}{M} \frac{1}{(N/M)}$$
$$\simeq \frac{(1 - q_k)\lambda_0(A)}{\sum_{r=1}^{K} p_r(1 - q_r)}.$$

Let
$$c = \max_{1 \leq k \leq K} \frac{(1 - q_k)}{\sum_{r=1}^{K} p_r(1 - q_r)}.$$

Then, for each $i \in I$ and any $A \in \mathcal{I}_0$, $P(\pi_i^{-1}(A)) \leq c \cdot \lambda(A)$. Thus, Keisler's Fubini property as in [26] and [34] also implies that $(\lambda \boxtimes P)(\pi^{-1}(A)) \leq c \cdot$ $\lambda(A)$. Let $B = \{i \in I : \alpha(i) \neq \beta(i)\}$. We know that $\lambda(B) = 0$, $(\lambda \boxtimes P)(\pi^{-1}(B)) = P(\pi_i^{-1}(B)) = 0$ for each $i \in I$. Since $g$ and $g^\beta$ agree on $I \times \Omega - \pi^{-1}(B)$, $g^\beta$ must be $\mathcal{I} \boxtimes \mathcal{F}$-measurable. For each $i \in I$, $g_i$ and $g_i^\beta$ agree on $\Omega - \pi_i^{-1}(B)$, and hence the relevant distribution and independence conditions are also satisfied by $g^\beta$.

REMARK 4.1. The sample space $\Omega$ in the proof of Theorem 2.6 depends on the choice of the internal type function $\alpha$. In the proof of Theorem 3.1 in Section 6 below, it will be more convenient to construct a sample space $\Omega$ that depends only on the agent space $I$, and not on the type function $\alpha$.

Let $\bar{\Omega}$ be the set of all internal bijections $\sigma$ from $I$ to $I$ such that for each $i \in I$, $\sigma(i) = i$ or $\sigma(\sigma(i)) = i$, and let $\bar{\mathcal{F}}_0$ be its internal power set. Let $q_k \in [0,1]$ for each $1 \leq k \leq K$. We can define an internal one-to-one mapping $\varphi^\alpha$ from $\Omega$ to $\bar{\Omega}$ by letting $\varphi^\alpha(B_1, B_2, \ldots, B_K, \pi^{B_1, B_2, \ldots, B_K})$ be the internal bijection $\sigma$ on $I$ such that $\sigma(i) = i$ for $i \in \bigcup_{k=1}^K B_k$ and $\sigma(i) = \pi^{B_1, B_2, \ldots, B_K}(i)$ for $i \in I - \bigcup_{k=1}^K B_k$. The mapping $\varphi^\alpha$ also induces an internal probability measure $\bar{P}_0^\alpha$ on $(\bar{\Omega}, \bar{\mathcal{F}}_0)$.

Define a mapping $\bar{\pi} : (I \times \bar{\Omega}) \to (I \cup J)$ by letting $\bar{\pi}(i, \bar{\omega}) = J$ if $\bar{\omega}(i) = i$, and $\bar{\pi}(i, \bar{\omega}) = \bar{\omega}(i)$ if $\bar{\omega}(i) \neq i$. Extend $\alpha$ so that $\alpha(J) = J$ and define $\bar{g}^\alpha : (I \times \bar{\Omega}) \to (S \cup \{J\})$ by letting

$$\bar{g}^\alpha(i, \bar{\omega}) = \alpha(\bar{\pi}(i, \bar{\omega}), \bar{\omega}).$$

Then it is obvious that $\bar{\pi}$ is still an independent-in-types random partial matching $\pi$ from $(I \times \bar{\Omega}, \mathcal{I} \boxtimes \bar{\mathcal{F}}, \lambda \boxtimes \bar{P}^\alpha)$ to $I$ with $q = (q_1, \ldots, q_K)$ as the no-match probabilities.

**5. Generalized Fubini and Ionescu–Tulcea theorems for Loeb transition probabilities.** For the proof of Theorem 3.1 to follow in Section 6, we need to work with a Loeb product transition probability system for a sequence of internal transition probabilities, based on the Loeb space construction in [33]. We first prove a generalized Fubini theorem for a Loeb transition probability in Section 5.1. Then, a generalized Ionescu–Tulcea theorem for an infinite sequence of Loeb transition probabilities is shown in Section 5.2. Finally, a Fubini extension based on Loeb product transition probability system is constructed in Section 5.3.

5.1. *A generalized Fubini theorem for a Loeb transition probability.* Let $(I, \mathcal{I}_0, \lambda_0)$ be a hyperfinite internal probability space with $\mathcal{I}_0$ the internal power set on a hyperfinite set $I$, and $\Omega$ a hyperfinite internal set with $\mathcal{F}_0$ its internal power set. Let $P_0$ be an internal function from $I$ to the space of hyperfinite internal probability measures on $(\Omega, \mathcal{F}_0)$, which is called an internal transition probability. For $i \in I$, denote the hyperfinite internal probability measure $P_0(i)$ by $P_{0i}$. Define a hyperfinite internal probability measure $\tau_0$ on $(I \times \Omega, \mathcal{I}_0 \otimes \mathcal{F}_0)$ by letting $\tau_0(\{(i,\omega)\}) = \lambda_0(\{i\}) P_{0i}(\{\omega\})$ for $(i,\omega) \in I \times \Omega$. Let $(I, \mathcal{I}, \lambda)$, $(\Omega, \mathcal{F}_i, P_i)$ and $(I \times \Omega, \mathcal{I} \boxtimes \mathcal{F}, \tau)$ be the Loeb spaces corresponding respectively to $(I, \mathcal{I}_0, \lambda_0)$, $(\Omega, \mathcal{F}_0, P_{0i})$ and $(I \times \Omega, \mathcal{I}_0 \otimes \mathcal{F}_0, \tau_0)$. The collection $\{P_i : i \in I\}$ of Loeb measures will be called a Loeb transition probability, and denoted by $P$. The measure $\tau$ will be called the Loeb product of the measure $\lambda$ and the Loeb transition probability $P$. We shall also denote $\tau_0$ by $\lambda_0 \otimes P_0$ and $\tau$ by $\lambda \boxtimes P$.

When $P_{0i} = P_0'$ for some hyperfinite internal probability measure $P_0'$ on $(\Omega, \mathcal{F}_0)$, $\tau_0$ is simply the internal product of $\lambda_0$ and $P_0'$, and $\tau$ the Loeb product measure $\lambda \boxtimes P'$, where $P'$ is the Loeb measure of $P_0'$. A Fubini-type theorem for this special case was shown by Keisler in [26], which is often referred to as Keisler's Fubini theorem. The following theorem presents a Fubini-type theorem for the general case.

THEOREM 5.1. *Let $f$ be a real-valued integrable function on $(I \times \Omega, \sigma(\mathcal{I}_0 \otimes \mathcal{F}_0), \tau)$. Then, (1) $f_i$ is $\sigma(\mathcal{F}_0)$-measurable for each $i \in I$ and integrable on $(\Omega, \sigma(\mathcal{F}_0), P_i)$ for $\lambda$-almost all $i \in I$; (2) $\int_\Omega f_i(\omega) \, dP_i(\omega)$ is integrable on $(I, \sigma(\mathcal{I}_0), \lambda)$; (3) $\int_I \int_\Omega f_i(\omega) \, dP_i(\omega) \, d\lambda(i) = \int_{I \times \Omega} f(i,\omega) \, d\tau(i,\omega)$.*

PROOF. Let $\mathcal{H}$ be the class of functions $g$ from $I \times \Omega$ to $\mathbb{R}_+ \cup \{+\infty\}$ that satisfy (1) for every $i \in I$, $g_i(\cdot)$ is $\sigma(\mathcal{F}_0)$-measurable; (2) the integral $\int_\Omega g_i(\omega) \, dP_i(\omega)$ as a function from $I$ to $\mathbb{R}_+ \cup \{+\infty\}$ is $\sigma(\mathcal{I}_0)$-measurable; (3) $\int_I \int_\Omega g_i(\omega) \, dP_i(\omega) \, d\lambda(i) = \int_{I \times \Omega} g(i,\omega) \, d\tau(i,\omega)$. It is obvious that $\mathcal{H}$ is closed under nonnegative linear combinations and monotone convergence.

Now, we consider $E \in \mathcal{I}_0 \otimes \mathcal{F}_0$ and $g = 1_E$. Then, for each $i \in I$, $g_i$ is the indicator function of the internal set $E_i = \{\omega \in \Omega : (i,\omega) \in E\}$, which is $\mathcal{F}_0$-measurable [and hence $\sigma(\mathcal{F}_0)$-measurable]. The integral $\int_\Omega g_i(\omega) \, dP_i(\omega) = P_i(E_i)$ is the standard part $°(P_{0i}(E_i))$. Since $P_{0i}(E_i)$ is $\mathcal{I}_0$-measurable as a function on $I$, $P_i(E_i)$ is thus $\sigma(\mathcal{I}_0)$-measurable as a function on $I$. Thus, the usual result on $S$-integrability (see, e.g., [34], Theorem 5.3.5, page 155) implies that

$$\text{(32)} \qquad \int_I \int_\Omega g_i(\omega) \, dP_i(\omega) \, d\lambda(i) = \int_I {}° P_{0i}(E_i) \, d\lambda(i)$$

$$\text{(33)} \qquad\qquad\qquad\qquad = {}°\!\int_I P_{0i}(E_i) \, d\lambda_0(i)$$

$$\text{(34)} \qquad\qquad\qquad\qquad = {}°\tau_0(E) = \int_{I \times \Omega} g \, d\tau.$$

Thus, $g \in \mathcal{H}$, and hence $\mathcal{H}$ contains the algebra $\mathcal{I}_0 \otimes \mathcal{F}_0$.

Therefore $\mathcal{H}$ is a monotone class. Then Theorem 3 on page 16 of [7] implies that $\mathcal{H}$ must contain all the nonnegative $\sigma(\mathcal{I}_0 \otimes \mathcal{F}_0)$-measurable functions.[22]

Since $f$ is integrable on $(I \times \Omega, \sigma(\mathcal{I}_0 \otimes \mathcal{F}_0), \tau)$, so are both $f^+$ and $f^-$. Now the fact that $\int_I \int_\Omega f_i^+(\omega) \, dP_i(\omega) \, d\lambda(i) < \infty$ implies that for $\lambda$-almost all $i \in I$, $\int_\Omega f_i^+(\omega) \, dP_i(\omega) < \infty$, and thus the $\sigma(\mathcal{F}_0)$-measurable function $f_i^+$ is integrable. Similarly, the $\sigma(\mathcal{I}_0)$-measurable function $\int_\Omega f_i^+(\omega) \, dP_i(\omega)$ is integrable since $\int_I \int_\Omega f_i^+(\omega) \, dP_i(\omega) \, d\lambda(i) < \infty$. We have similar results for $f^-$. The rest is clear. □

For any $B \in \sigma(\mathcal{F}_0)$, apply Theorem 5.1 to $f = 1_{I \times B}$ to obtain that $P_i(B)$ is $\sigma(\mathcal{I}_0)$-measurable for each $i \in I$. This means that $P = \{P_i : i \in I\}$ is indeed a transition probability in the usual sense (see [25]). One can define its usual product $\lambda \otimes P$ with $\lambda$ by letting $\lambda \otimes P(E) = \int_I P_i(E_i) \, d\lambda(i)$ for each $E$ in the usual product $\sigma$-algebra $\sigma(\mathcal{I}_0) \otimes \sigma(\mathcal{F}_0)$. It is clear that $(I \times \Omega, \sigma(\mathcal{I}_0 \otimes \mathcal{F}_0), \lambda \boxtimes P)$ is an extension of $(I \times \Omega, \sigma(\mathcal{I}_0) \otimes \sigma(\mathcal{F}_0), \lambda \otimes P)$.

The following result extends Theorem 5.1 to integrable functions on $(I \times \Omega, \mathcal{I} \boxtimes \mathcal{F}, \tau)$, which is the completion of $(I \times \Omega, \sigma(\mathcal{I}_0 \otimes \mathcal{F}_0), \tau)$.

PROPOSITION 5.2. *Let $f$ be a real-valued integrable function on $(I \times \Omega, \mathcal{I} \boxtimes \mathcal{F}, \tau)$. Then, (1) for $\lambda$-almost all $i \in I$, $f_i$ is integrable on the Loeb space $(\Omega, \mathcal{F}_i, P_i)$; (2) $\int_\Omega f_i(\omega) \, dP_i(\omega)$ is integrable on $(I, \mathcal{I}, \lambda)$; and (3) we have*

$$\int_I \int_\Omega f_i(\omega) \, dP_i(\omega) \, d\lambda(i) = \int_{I \times \Omega} f(i, \omega) \, d\tau(i, \omega).$$

PROOF. First, let $E \in \mathcal{I} \boxtimes \mathcal{F}$ with $\tau(E) = 0$. Then there is a set $A \in \sigma(\mathcal{I}_0 \otimes \mathcal{F}_0)$ such that $E \subseteq A$ and $\tau(A) = 0$. By Theorem 5.1, for $\lambda$-almost all $i \in I$, $P_i(A_i) = 0$, which implies that $P_i(E_i) = 0$.

There is a real-valued $\sigma(\mathcal{I}_0 \otimes \mathcal{F}_0)$-measurable function $g$ on $I \times \Omega$ such that the set $E = \{(i, \omega) \in I \times \Omega : f(i, \omega) \neq g(i, \omega)\}$ has $\tau$-measure zero. The above result implies that for $\lambda$-almost all $i \in I$, $f_i$ is the same as $g_i$. The rest is clear. □

Let $(I', \mathcal{I}'_0, \lambda'_0)$ be a hyperfinite internal probability space with $\mathcal{I}'_0$ the internal power set on a hyperfinite set $I'$. One can define a new transition probability, $\{\lambda'_0 \otimes P_{0i} : i \in I\}$. Define a hyperfinite internal probability measure $\tau_0^1$ on $(I' \times I \times \Omega, \mathcal{I}'_0 \otimes \mathcal{I}_0 \otimes \mathcal{F}_0)$ by letting

$$\tau_0^1(\{(i', i, \omega)\}) = \lambda_0(\{i\})(\lambda'_0 \otimes P_{0i})(\{(i', \omega)\}) = \lambda_0(\{i\}) \lambda'_0(\{i'\}) P_{0i}(\{\omega\})$$

for $(i', i, \omega) \in I' \times I \times \Omega$. Then it is clear that $\tau_0^1$ is exactly the same as $\lambda'_0 \otimes \tau_0$. Let $\lambda' \boxtimes \tau$ be the Loeb measure of $\lambda'_0 \otimes \tau_0$. By applying Theorem 5.1, we can obtain the following corollary.

---

[22]There is a typo in Theorem 3 on page 16 of [7]; $\mathcal{D}$ should be an algebra (not a $\sigma$-algebra as stated).

COROLLARY 5.3. *Let $h$ be a real-valued integrable function on $(I' \times I \times \Omega, \sigma(\mathcal{I}'_0 \otimes \mathcal{I}_0 \otimes \mathcal{F}_0), \lambda' \boxtimes \tau)$. Then, the following results hold.*

(1) $h_i$ *is $\sigma(\mathcal{I}'_0 \otimes \mathcal{F}_0)$-measurable for each $i \in I$ and integrable on $(\Omega, \sigma(\mathcal{I}'_0 \otimes \mathcal{F}_0), \lambda' \boxtimes P_i)$ for $\lambda$-almost all $i \in I$, where $\lambda' \boxtimes P_i$ is the Loeb measure of $\lambda'_0 \otimes P_{0i}$.*

(2) $\int_{I' \times \Omega} h_i(i', \omega) \, d\lambda' \boxtimes P_i(i', \omega)$ *is integrable on $(I, \sigma(\mathcal{I}_0), \lambda)$.*

(3) $\int_I \int_{I' \times \Omega} h_i(i', \omega) \, d(\lambda' \boxtimes P_i)(i', \omega) \, d\lambda(i) = \int_{I' \times I \times \Omega} h(i', i, \omega) \, d(\lambda' \boxtimes \tau)(i', i, \omega)$.

One can also view $\{P_{0i} : i \in I\}$ as a transition probability from $I' \times I$ to $\Omega$. Then the following corollary is obvious.

COROLLARY 5.4. *Let $h$ be a real-valued integrable function on $(I' \times I \times \Omega, \sigma(\mathcal{I}'_0 \otimes \mathcal{I}_0 \otimes \mathcal{F}_0), \lambda' \boxtimes \tau)$. Then, the following results hold.*

(1) $h_{(i',i)}$ *is $\sigma(\mathcal{F}_0)$-measurable for each $(i', i) \in I' \times I$ and integrable on $(\Omega, \sigma(\mathcal{F}_0), P_i)$ for $\lambda' \boxtimes \lambda$-almost all $(i', i) \in I' \times I$, where $\lambda' \boxtimes \lambda$ is the Loeb measure of $\lambda'_0 \otimes \lambda_0$.*

(2) $\int_\Omega h_{(i',i)}(\omega) \, dP_i(\omega)$ *is integrable on $(I' \times I, \sigma(\mathcal{I}'_0 \otimes \mathcal{I}_0), \lambda' \boxtimes \lambda)$.*

(3) $\int_{I' \times I} \int_\Omega h_{(i',i)}(\omega) \, dP_i(\omega) \, d(\lambda' \boxtimes \lambda)(i', i) = \int_{I' \times I \times \Omega} h(i', i, \omega) \, d(\lambda' \boxtimes \tau)(i', i, \omega)$.

By applying Proposition 5.2, one can also extend the results in Corollaries 5.3 and 5.4 to integrable functions on the Loeb space of $(I' \times I \times \Omega, \mathcal{I}'_0 \otimes \mathcal{I}_0 \otimes \mathcal{F}_0, \tau_0^1)$.

5.2. *A generalized Ionescu–Tulcea theorem for a Loeb product transition probability system.* In Section 5.1, the subscript 0 is used to distinguish internal measures and algebras with their corresponding Loeb measures and Loeb algebras. In this subsection, we need to work with infinitely many internal measure spaces and the corresponding Loeb spaces. To avoid confusion with the notation, we shall use $Q$ with subscripts or superscripts to represent internal measures; when their corresponding Loeb measures are considered, we use $P$ to replace $Q$. Internal algebras are denoted by $\mathcal{F}$ with subscripts or superscripts, and their external versions by $\mathcal{A}$ with subscripts or superscripts.

For each $m \geq 1$, let $\Omega_m$ be a hyperfinite set with its internal power set $\mathcal{F}_m$. We shall use $\Omega^n$, $\Omega^\infty$ and $\Omega_n^\infty$ to denote $\prod_{m=1}^n \Omega_m$, $\prod_{m=1}^\infty \Omega_m$ and $\prod_{m=n}^\infty \Omega_m$, respectively; also $\{\omega_m\}_{m=1}^n$, $\{\omega_m\}_{m=1}^\infty$ and $\{\omega_m\}_{m=n}^\infty$ will be denoted respectively by $\omega^n$, $\omega^\infty$ and $\omega_n^\infty$ when there is no confusion.

For each $n \geq 1$, let $Q_n$ be an internal transition probability from $\Omega^{n-1}$ to $(\Omega_n, \mathcal{F}_n)$; that is, for each $\omega^{n-1} \in \Omega^{n-1}$, $Q_n(\omega^{n-1})$ (also denoted by $Q_n^{\omega^{n-1}}$) is a hyperfinite internal probability measure on $(\Omega_n, \mathcal{F}_n)$. In particular, $Q_1$ is simply a hyperfinite internal probability measure on $(\Omega_1, \mathcal{F}_1)$, and $Q_2$ an internal transition probability from $\Omega_1$ to $(\Omega_2, \mathcal{F}_2)$. Thus, $Q_1 \otimes Q_2$ defines an internal probability measure on $(\Omega_1 \times \Omega_2, \mathcal{F}_1 \otimes \mathcal{F}_2)$. By induction, $Q_1 \otimes Q_2 \otimes \cdots \otimes Q_n$ defines an internal probability measure on $(\Omega^n, \bigotimes_{m=1}^n \mathcal{F}_m)$.[23] Denote $Q_1 \otimes Q_2 \otimes \cdots \otimes Q_n$ by $Q^n$, and $\bigotimes_{m=1}^n \mathcal{F}_m$ by $\mathcal{F}^n$. Then $Q^n$ is the internal product of the internal transition probability $Q_n$ with the internal probability measure $Q^{n-1}$. Let $P^n$ and $P_n(\omega^{n-1})$ (also denoted by $P_n^{\omega^{n-1}}$) be the corresponding Loeb measures, which are defined respectively on $\sigma(\mathcal{F}^n)$ and $\sigma(\mathcal{F}_n)$. Using the notation in Section 5.1, $P^n$ is the Loeb product $P_1 \boxtimes P_2 \boxtimes \cdots \boxtimes P_n$ of the Loeb transition probabilities $P_1, P_2, \ldots, P_n$.

Theorem 5.1 implies that for any set $E \in \sigma(\mathcal{F}^n)$,

$$(35) \qquad P^n(E) = \int_{\Omega^{n-1}} P_n^{\omega^{n-1}}(E_{\omega^{n-1}}) \, dP^{n-1}(\omega^{n-1}).$$

That is, $P^n$ is the product of the transition probability $P_n$ with the probability measure $P^{n-1}$. Thus, we can also denote $P^n$ by $P^{n-1} \boxtimes P_n$, and furthermore by $\boxtimes_{m=1}^n P_m$.

Let $\pi_n^{n-1}$ be the projection mapping from $\Omega^n$ to $\Omega^{n-1}$; that is, $\pi_n^{n-1}(\omega_1, \ldots, \omega_n) = (\omega_1, \ldots, \omega_{n-1})$. Let $F$ be any subset of $\Omega^{n-1}$ and $E = F \times \Omega_n$. Then, $E_{\omega^{n-1}} = \Omega_n$ when $\omega^{n-1} \in F$, and $E_{\omega^{n-1}} = \varnothing$ when $\omega^{n-1} \notin F$. If $E \in \sigma(\mathcal{F}^n)$, then Theorem 5.1 implies that $P_n^{\omega^{n-1}}(E_{\omega^{n-1}}) = 1_F$ is $\sigma(\mathcal{F}^{n-1})$-measurable [i.e., $F \in \sigma(\mathcal{F}^{n-1})$], and $P^n(E) = P^{n-1}(F)$. On the other hand, if $F \in \sigma(\mathcal{F}^{n-1})$, then it is obvious that $E \in \sigma(\mathcal{F}^n)$ and $P^n(E) = P^{n-1}(F)$. This means that the measure space $(\Omega^{n-1}, \sigma(\mathcal{F}^{n-1}), P^{n-1})$ is the projection of $(\Omega^n, \sigma(\mathcal{F}^n), P^n)$ under $\pi_n^{n-1}$. Similarly, let $\pi_n^k$ be the projection mapping from $\Omega^n$ to $\Omega^k$ for some $k < n$; then $(\Omega^k, \sigma(\mathcal{F}^k), P^k)$ is the projection of $(\Omega^n, \sigma(\mathcal{F}^n), P^n)$ under $\pi_n^k$.

For a collection $\mathcal{D}$ of sets and a set $F$, we use $\mathcal{D} \times F$ to denote the collection $\{D \times F : D \in \mathcal{D}\}$ of sets when there is no confusion. Thus, $\sigma(\mathcal{F}^n) \times \Omega_{n+1}^\infty$ denotes $\{E_n \times \Omega_{n+1}^\infty : E_n \in \sigma(\mathcal{F}^n)\}$. Let $\mathcal{E} = \bigcup_{n=1}^\infty [\sigma(\mathcal{F}^n) \times \Omega_{n+1}^\infty]$, which is an algebra of sets in $\Omega^\infty$. One can define a measure $P^\infty$ on this algebra by letting $P^\infty(E_n \times \Omega_{n+1}^\infty) = P^n(E_n)$ for each $E_n \in \sigma(\mathcal{F}^n)$. The projection property stated in the above paragraph implies that $P^\infty$ is well defined. Let $\mathcal{F}^\infty = \bigcup_{n=1}^\infty [\mathcal{F}^n \times \Omega_{n+1}^\infty]$. Then, it is clear that $\sigma(\mathcal{F}^\infty) = \sigma(\mathcal{E})$.

The point is how to extend $P^\infty$ to a countably additive probability measure on the $\sigma$-algebra $\sigma(\mathcal{F}^\infty)$. This is possible by using a proof similar to that

---

[23]In fact, for each $\omega^n \in \Omega^n$, we have

$$\bigotimes_{m=1}^n Q_m(\{(\omega_1, \ldots, \omega_n)\}) = \prod_{m=1}^n Q_m(\omega_1, \ldots, \omega_{m-1})(\{\omega_m\}).$$

of Proposition 3.3 of [39]. The result is a version of the Ionescu–Tulcea theorem (see [25], page 93) for the Loeb product transition probability system $\{P_1 \boxtimes P_2 \boxtimes \cdots \boxtimes P_n\}_{n=1}^\infty$.

THEOREM 5.5. *There is a unique countably additive probability measure on $\sigma(\mathcal{F}^\infty)$ that extends the set function $P^\infty$ on $\mathcal{E}$; such a unique extension is still denoted by $P^\infty$ and by $\boxtimes_{m=1}^\infty P_m$.*

PROOF. Let $\{C_n\}_{n=1}^\infty$ be a decreasing sequence of sets in $\mathcal{F}^\infty$ with empty intersection. By the construction of $\mathcal{F}^\infty$, one can find a sequence of internal sets $\{A_n\}_{n=1}^\infty$ and a nondecreasing sequence $\{k_n\}_{n=1}^\infty$ of nonnegative integers such that $C_n = A_n \times \Omega_{k_n+1}^\infty$ and $A_n \in \mathcal{F}^{k_n}$. For $\ell \leq n$, let $\pi_{k_n}^{k_l}$ be the mapping from $\Omega^{k_n}$ to $\Omega^{k_l}$ by projecting a tuple in $\Omega^{k_n}$ to its first $k_\ell$ coordinates; then $\pi_{k_n}^{k_l}(A_n) \subseteq A_\ell$ because $\{C_n\}_{n=1}^\infty$ is a decreasing sequence of sets. Take the transfer $\{k_n\}_{n\in{}^*\mathbb{N}}$ of the sequence $\{k_n\}_{n=1}^\infty$, and the respective internal extensions $\{A_n\}_{n\in{}^*\mathbb{N}}$ and $\{\Omega^n\}_{n\in{}^*\mathbb{N}}$ of the internal sequences $\{A_n\}_{n=1}^\infty$ and $\{\Omega^n\}_{n\in\mathbb{N}}$. By spillover and $\aleph_1$-saturation (see [34]), one can obtain $h \in {}^*\mathbb{N}_\infty$ such that for all $n \leq h$, $A_n \subseteq \Omega^{k_n}$ and $\pi_{k_n}^{k_l}(A_n) \subseteq A_\ell$ for all $\ell \in {}^*\mathbb{N}$ with $l \leq n$, where $\pi_{k_n}^{k_l}$ is defined in exactly the same way as in the case of finite $n$.

We claim that $A_n = \varnothing$ for all $n \in {}^*\mathbb{N}_\infty$ with $n \leq h$; if not, one can find such an $n$ with $\omega^{k_n} = (\omega_1, \ldots, \omega_n) \in A_n$. Then $\omega^{k_l} \in A_l$ for any $\ell \in \mathbb{N}$. If $k_n \in {}^*\mathbb{N}_\infty$, then it is obvious that $\{\omega_m\}_{m=1}^\infty$ is in $C_\ell$ for all $\ell \in \mathbb{N}$, which contradicts the assumption that the intersection of all the $C_\ell$ is empty. If $k_n \in \mathbb{N}$, one can choose $\omega_m$ arbitrarily for any $m > k_n$ to obtain the same contradiction. Hence the claim is proven.

By spillover, we know that for some $n \in \mathbb{N}$, $A_n = \varnothing$, and so is $C_n$. Thus, we obtain a trivial limit, $\lim_{n\to\infty} P^\infty(C_n) = 0$. This means that $P^\infty$ is indeed countably additive on $\mathcal{F}^\infty$. As in [33], the Carathéodory extension theorem implies that $P^\infty$ can be extended to the $\sigma$-algebra $\sigma(\mathcal{F}^\infty)$ generated by $\mathcal{F}^\infty$, and we are done. □

The following result considers sectional measurability for sets in $\sigma(\mathcal{F}^\infty)$.

PROPOSITION 5.6. *Let $G$ be a $\sigma(\mathcal{F}^\infty)$-measurable subset of $\Omega^\infty$. Then, for any $\{\omega_m\}_{m=1}^\infty \in \Omega^\infty$, the set $G_{\omega_{n+1}^\infty} = \{\omega'^n \in \Omega^n : (\omega_1', \ldots, \omega_n', \omega_{n+1}, \omega_{n+2}, \ldots) \in G\}$ belongs to $\sigma(\mathcal{F}^n)$, while the set*

$$G_{\omega^n} = \{\omega'^\infty_{n+1} \in \Omega_{n+1}^\infty : (\omega_1, \ldots, \omega_n, \omega'_{n+1}, \omega'_{n+2}, \ldots) \in G\}$$

*belongs to $\sigma(\bigcup_{m=n+1}^\infty [(\bigotimes_{k=n+1}^m \mathcal{F}_k) \times \Omega_{m+1}^\infty])$.*

PROOF. The collection of those sets $G$ in $\sigma(\mathcal{F}^\infty)$ with the properties is clearly a monotone class of sets and contains the algebra $\mathcal{F}^\infty$; and hence it is $\sigma(\mathcal{F}^\infty)$ itself by Theorem 1 on page 7 of [7]. □

The following corollary follows from Proposition 5.6 immediately.

COROLLARY 5.7. *Let $\pi^n$ be the projection mapping from $\Omega^\infty$ to $\Omega^n$ [i.e., $\pi^n(\omega^\infty) = \omega^n$]. Then the measure space $(\Omega^n, \sigma(\mathcal{F}^n), P^n)$ is the projection of the measure space $(\Omega^\infty, \sigma(\mathcal{F}^\infty), P^\infty)$ under $\pi^n$ in the sense that for any $F \subseteq \Omega^n$, $F \in \sigma(\mathcal{F}^n)$ if and only if $(\pi^n)^{-1}(F) = F \times \Omega_{n+1}^\infty \in \sigma(\mathcal{F}^\infty)$ with $P^n(F) = P^\infty((\pi^n)^{-1}(F))$.*

5.3. *Fubini extensions based on a Loeb product transition probability system.* Let $(I, \mathcal{I}_0, \lambda_0)$ be a hyperfinite internal probability space with $\mathcal{I}_0$ the internal power set on a hyperfinite set $I$. We follow the notation and construction in Section 5.2. Denote $(\Omega_0, \mathcal{F}_0, Q_0) = (I, \mathcal{I}_0, \lambda_0)$ and repeat the process of constructing a countably additive measure $\boxtimes_{m=0}^\infty P_m$ on $(I \times \Omega^\infty, \sigma(\bigcup_{n=1}^\infty (\mathcal{I}_0 \otimes \mathcal{F}^n) \times \Omega_{n+1}^\infty))$.

The following lemma is a restatement of Keisler's Fubini theorem to the particular setting. Since the marginals of $\boxtimes_{m=0}^n P_m$ on $I$ and $\Omega^n$ are respectively $\lambda$ and $P^n$, we can write $\boxtimes_{m=0}^n P_m$ as $\lambda \boxtimes P^n$.

LEMMA 5.8. *For any $n \geq 1$, the space $(I \times \Omega^n, \sigma(\mathcal{I}_0 \otimes \mathcal{F}^n), \boxtimes_{m=0}^n P_m)$ is a Fubini extension over the usual product of $(I, \sigma(\mathcal{I}_0), \lambda)$ and $(\Omega^n, \sigma(\mathcal{F}^n), P^n)$.*

The following is a Fubini-type result for the infinite product.

PROPOSITION 5.9. *The space $(I \times \Omega^\infty, \sigma(\bigcup_{n=1}^\infty (\mathcal{I}_0 \otimes \mathcal{F}^n) \times \Omega_{n+1}^\infty), \boxtimes_{m=0}^\infty P_m)$ is a Fubini extension over the usual product of the probability spaces $(I, \mathcal{I}, \lambda)$ and $(\Omega^\infty, \sigma(\mathcal{F}^\infty), P^\infty)$.*

PROOF. We only check that the Fubini-type property holds for any set $E \in \bigcup_{n=1}^\infty (\mathcal{I}_0 \otimes \mathcal{F}^n) \times \Omega_{n+1}^\infty$. The rest of the proof is essentially the same as that of Theorem 5.1.

It is clear that there exists a set $F \in \mathcal{I}_0 \otimes \mathcal{F}^n$ such that $E = F \times \Omega_{n+1}^\infty$. By definition, $\boxtimes_{m=0}^\infty P_m(E) = \boxtimes_{m=0}^n P_m(F)$. By Lemma 5.8,

$$\underset{m=0}{\overset{n}{\boxtimes}} P_m(F) = \lambda \boxtimes P^n(F) = \int_I P^n(F_i)\,d\lambda(i) = \int_{\Omega^n} \lambda(F_{\omega^n})\,dP^n(\omega^n).$$

On the other hand, $E_i = F_i \times \Omega_{n+1}^\infty$ and $E_{\omega^\infty} = F_{\omega^n}$. By the fact that $P^\infty(E_i) = P^n(F_i)$, the projection property in Corollary 5.7 implies that

$$\underset{m=0}{\overset{n}{\boxtimes}} P_m(F) = \int_I P^\infty(F_i)\,d\lambda(i) = \int_{\Omega^\infty} \lambda(F_{\omega^\infty})\,dP^\infty(\omega^\infty).$$

This means that the Fubini property does hold for sets in $\bigcup_{n=1}^\infty (\mathcal{I}_0 \otimes \mathcal{F}^n) \times \Omega_{n+1}^\infty$. □

For simplicity, we denote $\sigma(\mathcal{I}_0)$, $\sigma(\mathcal{F}_n)$, $\sigma(\mathcal{I}_0 \otimes \mathcal{F}_n)$, $\sigma(\mathcal{F}^n)$, $\sigma(\mathcal{I}_0 \otimes \mathcal{F}^n)$, $\sigma(\mathcal{F}^\infty)$ and $\sigma(\bigcup_{n=1}^\infty (\mathcal{I}_0 \otimes \mathcal{F}^n) \times \Omega_{n+1}^\infty)$, respectively by $\mathcal{I}$, $\mathcal{A}_n$, $\mathcal{I} \boxtimes \mathcal{A}_n$, $\mathcal{A}^n$, $\mathcal{I} \boxtimes \mathcal{A}^n$, $\mathcal{A}^\infty$ and $\mathcal{I} \boxtimes \mathcal{A}^\infty$.

We restate some of the above results using the new notation. Corollary 5.7 implies that $(\Omega^n, \mathcal{A}^n, P^n)$ and $(I \times \Omega^n, \mathcal{I} \boxtimes \mathcal{A}^n, \lambda \boxtimes P^n)$ are the respective projections of $(\Omega^\infty, \mathcal{A}^\infty, P^\infty)$ and $(I \times \Omega^\infty, \mathcal{I} \boxtimes \mathcal{A}^\infty, \lambda \boxtimes P^\infty)$. Since $(I \times \Omega^n, \mathcal{I} \boxtimes \mathcal{A}^n, \lambda \boxtimes P^n)$ is the Loeb product of two Loeb probability spaces $(I, \mathcal{I}, \lambda)$ and $(\Omega^n, \mathcal{A}^n, P^n)$, it is a Fubini extension of the usual product probability space. In addition, the Fubini property in Proposition 5.9 says that $(I \times \Omega^\infty, \mathcal{I} \boxtimes \mathcal{A}^\infty, \lambda \boxtimes P^\infty)$ is a Fubini extension of the usual product of the probability spaces $(I, \mathcal{I}, \lambda)$ and $(\Omega^\infty, \mathcal{A}^\infty, P^\infty)$.

**6. Proof of Theorem 3.1.** Let $(p^0, b, q, \nu)$ be the given parameters for the dynamical system $\mathbb{D}$. Let $M$ be a fixed unlimited hyperfinite natural number in $^*\mathbb{N}_\infty$, $I = \{1, 2, \ldots, M\}$, $\mathcal{I}_0$ the internal power set on $I$, and $\lambda_0$ the internal counting probability measure on $\mathcal{I}_0$. Let $\alpha^0 : I \to S = \{1, 2, \ldots, K\}$ be an internal initial type function such that $\lambda_0(\alpha^0 = k) \simeq p_k$ for each $k = 1, \ldots, K$.[24] What we need to do is to construct a sequence of internal transition probabilities and a sequence of internal type functions. The results in Sections 5.2 and 5.3 can then be applied to obtain a Loeb product transition probability system. Since we need to consider random mutation, random partial matching and random type changing at each time period, three internal measurable spaces with internal transition probabilities will be constructed at each time period.

Adopt the notation used in Section 5.2. Suppose that the construction for the dynamical system $\mathbb{D}$ has been done up to time period $n - 1$. Thus, $\{(\Omega_m, \mathcal{F}_m, Q_m)\}_{m=1}^{3n-3}$ and $\{\alpha^l\}_{l=0}^{n-1}$ have been constructed, where each $\Omega_m$ is a hyperfinite internal set with its internal power set $\mathcal{F}_m$, $Q_m$ an internal transition probability from $\Omega^{m-1}$ to $(\Omega_m, \mathcal{F}_m)$, and $\alpha^l$ an internal function from $I \times \Omega^{3l}$ to the type space $S$.

We shall now consider the constructions for time $n$. We first work with the random mutation step. For each $1 \le k \le K$, $\rho_k$ is a distribution on $S$ with $\rho_k(l) = b_{kl}$, the probability for a type-$k$ agent to mutate to a type-$l$ agent. Let $\Omega_{3n-2} = S^I$ (the space of all internal functions from $I$ to $S$) with its internal power set $\mathcal{F}_{3n-2}$.

For each $i \in I$, $\omega^{3n-3} \in \Omega^{3n-3}$, let $\gamma_i^{\omega^{3n-3}} = \rho_{\alpha^{n-1}(i, \omega^{3n-3})}$. That is, if $\alpha^{n-1}(i, \omega^{3n-3}) = k$, then $\gamma_i^{\omega^{3n-3}} = \rho_k$. Define an internal probability measure $Q_{3n-2}^{\omega^{3n-3}}$ on $(S^I, \mathcal{F}_{3n-2})$ to be the internal product measure $\prod_{i \in I} \gamma_i^{\omega^{3n-3}}$.[25] Let $h^n : (I \times \prod_{m=1}^{3n-2} \Omega_m) \to S$ be such that $h^n(i, \omega^{3n-2}) = \omega_{3n-2}(i)$.

---

[24]This is possible since $\lambda$ is atomless.

[25]A hyperfinite product space is a common construction in nonstandard analysis, whose coordinate functions also give a hyperfinite sequence of $*$-independent random variables.

Next, we consider the step of random partial matching. Let $(\Omega_{3n-1}, \mathcal{F}_{3n-1})$ be the internal sample measurable space $(\bar{\Omega}, \bar{\mathcal{F}}_0)$ in Remark 4.1.

For any given $\omega^{3n-2} \in \Omega^{3n-2}$, the type function is $h^n_{\omega^{3n-2}}(\cdot)$, denoted by $\alpha$ for short. Let $Q^{\omega^{3n-2}}_{3n-1}$ be the internal probability measure corresponding to the internal probability measure $\bar{P}^{\alpha}_0$ in Remark 4.1. Define a mapping $\pi^n : (I \times \Omega^{3n-1}) \to (I \cup J)$ by letting $\pi^n(i, \omega^{3n-1}) = \bar{\pi}(i, \omega_{3n-1})$; thus, $\pi^n(i, \omega^{3n-1}) = J$ if $\omega_{3n-1}(i) = i$, and $\pi^n(i, \omega^{3n-1}) = \omega_{3n-1}(i)$ if $\omega_{3n-1}(i) \neq i$. Extend $h^n$ so that $h^n(J, \omega^{3n-2}) = J$ for any $\omega^{3n-2} \in \Omega^{3n-2}$. Define $g^n : (I \times \Omega^{3n-1}) \to (S \cup \{J\})$ by letting

$$g^n(i, \omega^{3n-1}) = h^n(\pi^n(i, \omega^{3n-1}), \omega^{3n-2}),$$

which means that $g^n(i, \omega^{3n-1}) = \bar{g}^{\alpha}(i, \omega_{3n-1})$.

Finally, we consider the step of random type changing for matched agents. Let $\Omega_{3n} = S^I$ with its internal power set $\mathcal{F}_{3n}$; each point $\omega_{3n} \in \Omega_{3n}$ is an internal function from $I$ to $S$. For any given $\omega^{3n-1} \in \Omega^{3n-1}$, the space $I$ of agents is divided into $K^2 + K$ classes: those in type $k$ who are not matched, or matched to some type-$l$ agents. For $1 \leq k, l \leq K$, $\nu_{kl}$ is a distribution on $S$ and $\nu_{kl}(r)$ the probability for a type-$k$ agent to change to a type-$r$ agent when the type-$k$ agent meets a type-$l$ agent.

Define a new type function $\alpha^n : (I \times \Omega^{3n}) \to S$ by letting $\alpha^n(i, \omega^{3n}) = \omega_{3n}(i)$.

Fix $\omega^{3n-1} \in \Omega^{3n-1}$. For each $i \in I$, (1) if $\omega_{3n-1}(i) = i$ ($i$ is not matched at time $n$), let $\tau_i^{\omega^{3n-1}}$ be the probability measure on the type space $S$ that gives probability 1 to the type $h^n(i, \omega^{3n-2})$ and zero for the rest; (2) if $\omega_{3n-1}(i) \neq i$ ($i$ is matched at time $n$), $h^n(i, \omega^{3n-2}) = k$ and $h^n(\omega_{3n-1}(i), \omega^{3n-2}) = l$, let $\tau_i^{\omega^{3n-1}}$ be the distribution $\nu_{kl}$ on $S$. Define an internal probability measure $Q^{\omega^{3n-1}}_{3n}$ on $S^I$ to be the internal product measure $\prod_{i \in I} \tau_i^{\omega^{3n-1}}$.

By induction, we can construct a sequence $\{(\Omega_m, \mathcal{F}_m, Q_m)\}_{m=1}^{\infty}$ of internal transition probabilities and a sequence $\{\alpha^l\}_{l=0}^{\infty}$ of type functions. By using the constructions in Sections 5.2 and 5.3 via an infinite product of Loeb transition probabilities, we can obtain a corresponding probability space $(I \times \Omega^{\infty}, \mathcal{I} \boxtimes \mathcal{A}^{\infty}, \lambda \boxtimes P^{\infty})$.

From now on, we shall also use $(\Omega, \mathcal{F}, P)$ and $(I \times \Omega, \mathcal{I} \boxtimes \mathcal{F}, \lambda \boxtimes P)$ to denote $(\Omega^{\infty}, \mathcal{A}^{\infty}, P^{\infty})$ and $(I \times \Omega^{\infty}, \mathcal{I} \boxtimes \mathcal{A}^{\infty}, \lambda \boxtimes P^{\infty})$, respectively. Note that all the functions, $h^n, \pi^n, g^n, \alpha^n$, for $n = 1, 2, \ldots$, can be viewed as functions on $I \times \Omega$, and $h^n, g^n, \alpha^n$ are $\mathcal{I} \boxtimes \mathcal{F}$-measurable for each $n \geq 1$.

---

In general, for a hyperfinite collection $\{X_i\}_{i=1}^n$ (with $n$ infinitely large) of $*$-independent random variables on an internal probability space $(\Omega, \mathcal{F}_0, P_0)$ with mean zero and variances bounded by a common standard positive number $C$, the elementary Chebyshev's inequality says that for any positive hyperreal number $\varepsilon$, $P_0(|X_1 + \cdots + X_n|/n \geq \varepsilon) \leq C/n\varepsilon^2$. By taking $\varepsilon = 1/n^{1/4}$, it implies that for nearly all $\omega \in \Omega$, $|X_1 + \cdots + X_n|/n \simeq 0$, which was also noted in [26], page 56.

We still need to check that our internal constructions above lead to a dynamical system $\mathbb{D}$ with random mutation, partial matching and type changing that is Markov conditionally independent in types. We assume that the conditions for random mutation, partial matching and type changing as well as Markov conditional independence in types are satisfied up to time $n-1$. As in the proof of Lemma 6 in [15], Lemma 5 in [15] implies that the random variables $\alpha_i^{n-1}$ and $\alpha_j^{n-1}$ are independent for $i \neq j$. It remains to check the conditions for random mutation, partial matching and type changing as well as Markov conditional independence in types for time $n$.

For the step of random mutation at time period $n$, we have for each agent $i \in I$, and $k, l \in S$,

$$P(h_i^n = l, \alpha_i^{n-1} = k)$$
$$= P^{3n-2}(A_{3n-2})$$
(36)
$$= P^{3n-2}(\{\omega^{3n-3} \in \Omega^{3n-3} : \alpha_i^{n-1}(\omega^{3n-3}) = k\} \times (S^{I-\{i\}} \times \{l\}^{\{i\}}))$$
$$= \int_{\{\alpha_i^{n-1}(\omega^{3n-3}) = k\}} \rho_k(l) \, dP^{3n-3}(\omega^{3n-3}) = b_{kl} P(\alpha_i^{n-1} = k),$$

where

$$A_{3n-2} = \{(\omega^{3n-3}, \omega_{3n-2}) \in \Omega^{3n-2} : \alpha_i^{n-1}(\omega^{3n-3}) = k,$$
$$h^n(i, \omega^{3n-3}, \omega_{3n-2}) = \omega_{3n-2}(i) = l\},$$

which implies that (1) is satisfied.

When $i \neq j \in I$, it is obvious that for any $l, r \in S$, and any $(a_i^0, \ldots, a_i^{n-1})$ and $(a_j^0, \ldots, a_j^{n-1})$ in $S^n$, for the event

$$B_1 = \{h_i^n = l, h_j^n = r, (\alpha_i^0, \ldots, \alpha_i^{n-1}) = (a_i^0, \ldots, a_i^{n-1}),$$
$$(\alpha_j^0, \ldots, \alpha_j^{n-1}) = (a_j^0, \ldots, a_j^{n-1})\},$$

we have

$$P(B_1) = P^{3n-2}(B_2)$$
(37)
$$= P^{3n-2}(B_3 \times (I \times S^{I-\{i,j\}} \times \{l\}^{\{i\}} \times \{r\}^{\{j\}}))$$
$$= P^{3n-3}(B_4) \cdot \rho_{a_i^{n-1}}(l) \cdot \rho_{a_j^{n-1}}(r),$$

where

$$B_2 = \{(\omega^{3n-3}, \omega_{3n-2}) \in \Omega^{3n-2} : (\alpha_i^0, \ldots, \alpha_i^{n-1})(\omega^{3n-3}) = (a_i^0, \ldots, a_i^{n-1}),$$
$$(\alpha_j^0, \ldots, \alpha_j^{n-1})(\omega^{3n-3}) = (a_j^0, \ldots, a_j^{n-1}),$$
$$\omega_{3n-2}(i) = l, \omega_{3n-2}(j) = r\},$$
$$B_3 = \{\omega^{3n-3} \in \Omega^{3n-3} : (\alpha_i^0, \ldots, \alpha_i^{n-1})(\omega^{3n-3}) = (a_i^0, \ldots, a_i^{n-1}),$$
$$(\alpha_j^0, \ldots, \alpha_j^{n-1})(\omega^{3n-3}) = (a_j^0, \ldots, a_j^{n-1})\},$$

and

$$B_4 = \{\omega^{3n-3} \in \Omega^{3n-3} : (\alpha_m^0, \ldots, \alpha_m^{n-1})(\omega^{3n-3}) = (a_m^0, \ldots, a_m^{n-1}), m = i, j\}.$$

Equations (36) and (37) imply that for any $l, r \in S$,

$$P(h_i^n = l, h_j^n = r | (\alpha_m^0, \ldots, \alpha_m^{n-1}), m = i, j)$$
$$= P(h_i^n = l | \alpha_i^{n-1}) P(h_j^n = r | \alpha_j^{n-1}).$$

Hence (5) in the definition of Markov conditional independence for random mutation is satisfied.

Equation (37) together with the independence of $\alpha_i^{n-1}$ and $\alpha_j^{n-1}$ implies that $h_i^n(\cdot)$ and $h_j^n(\cdot)$ are independent. As in (2), $\overline{p}^{n-1/2}$ is the expected cross-sectional type distribution immediately after random mutation. The exact law of large numbers in [38] and [40] (see Footnote 20 above) implies that for $P^{3n-2}$-almost all $\omega^{3n-2} \in \Omega^{3n-2}$,

$$(38) \qquad \lambda(\{i' \in I : h_{\omega^{3n-2}}^n(i') = l\}) = \overline{p}_l^{n-1/2}$$

for any $l \in S$.

For the step of random partial matching, the definition of $\pi^n$ clearly shows that for each $\omega \in \Omega$, the restriction of $\pi_\omega^n(\cdot)$ to $I - (\pi_\omega^n)^{-1}(\{J\})$ is a full matching on that set. We need to check that the function $g^n : (I \times \Omega) \to (I \cup \{J\})$ satisfies the required distribution and Markov conditional independence conditions.

For any given $\omega^{3n-2} \in \Omega^{3n-2}$, take $\alpha = h_{\omega^{3n-2}}^n(\cdot)$ as in the construction $Q_{3n-1}^{\omega^{3n-2}}$ above. Then, for each $i \in I$ and $c \in S \cup \{J\}$,

$$(39) \qquad Q_{3n-1}^{\omega^{3n-2}}(\{\omega_{3n-1} \in \Omega_{3n-1} : g^n(i, \omega^{3n-1}) = c\})$$
$$= \bar{P}_0^\alpha(\{\bar{\omega} \in \bar{\Omega} : \bar{g}^\alpha(i, \bar{\omega}) = c\}).$$

Moreover, for each $j \in I$ with $j \neq i$ and $d \in S \cup \{J\}$,

$$Q_{3n-1}^{\omega^{3n-2}}(\{\omega_{3n-1} \in \Omega_{3n-1} : g^n(i, \omega^{3n-1}) = c, g^n(j, \omega^{3n-1}) = d\})$$
$$(40) \qquad = \bar{P}_0^\alpha(\{\bar{\omega} \in \bar{\Omega} : \bar{g}^\alpha(i, \bar{\omega}) = c, \bar{g}^\alpha(j, \bar{\omega}) = d\})$$
$$\simeq \bar{P}_0^\alpha(\bar{g}_i^\alpha = c) \cdot \bar{P}_0^\alpha(\bar{g}_j^\alpha = d).$$

For each agent $i \in I$, $k, l \in S$, (39) implies that

$$P(g_i^n = J, h_i^n = k)$$
$$= P^{3n-1}(\{\omega^{3n-1} \in \Omega^{3n-1} : h^n(i, \omega^{3n-2}) = k, g^n(i, \omega^{3n-1}) = J\})$$
$$(41) \qquad = \int_{\{h^n(i, \omega^{3n-2}) = k\}} q_k \, dP^{3n-2}(\omega^{3n-2})$$
$$= q_k P(h_i^n = k).$$

In addition, we obtain from (38) that

$$P(g_i^n = l, h_i^n = k)$$
$$= P^{3n-1}(\{\omega^{3n-1} \in \Omega^{3n-1} : h^n(i, \omega^{3n-2}) = k, g^n(i, \omega^{3n-1}) = l\})$$
$$= \int_{\{h^n(i,\omega^{3n-2})=k\}} \bar{P}^{h^n_{\omega^{3n-2}}}(\{\bar{\omega} \in \bar{\Omega} : \bar{g}^{h^n_{\omega^{3n-2}}}(i, \bar{\omega}) = l\}) \, dP^{3n-2}(\omega^{3n-2})$$

$$(42) = \int_{\{h^n(i,\omega^{3n-2})=k\}} \frac{(1-q_k)(1-q_l)\lambda(\{i' \in I : h^n_{\omega^{3n-2}}(i') = l\})}{\sum_{r=1}^{K}(1-q_r)\lambda(\{i' \in I : h^n_{\omega^{3n-2}}(i') = r\})} \, dP^{3n-2}(\omega^{3n-2})$$

$$= \int_{\{h^n(i,\omega^{3n-2})=k\}} \frac{(1-q_k)(1-q_l)\bar{p}_l^{n-1/2}}{\sum_{r=1}^{K}(1-q_r)\bar{p}_r^{n-1/2}} \, dP^{3n-2}(\omega^{3n-2})$$

$$= P(h_i^n = k) \cdot \frac{(1-q_k)(1-q_l)\bar{p}_l^{n-1/2}}{\sum_{r=1}^{K}(1-q_r)\bar{p}_r^{n-1/2}}.$$

Hence, (41) and (42) imply that (3) holds.

Fix $i \neq j \in I$. Take any $(a_i^0, \ldots, a_i^{n-1}), (a_j^0, \ldots, a_j^{n-1}) \in S^n$, $l, r \in S$, and $c, d \in S \cup \{J\}$. Let $D$ be the set of all $\omega^{3n-2} \in \Omega^{3n-2}$ such that $h_i^n(\omega^{3n-2}) = l$, $h_j^n(\omega^{3n-2}) = r$, $(\alpha_i^0, \ldots, \alpha_i^{n-1})(\omega^{3n-3}) = (a_i^0, \ldots, a_i^{n-1})$ and $(\alpha_j^0, \ldots, \alpha_j^{n-1})(\omega^{3n-3}) = (a_j^0, \ldots, a_j^{n-1})$. Then, (38) and (40) imply that

$$P(g_i^n = c, g_j^n = d, h_i^n = l, h_j^n = r,$$
$$(\alpha_m^0, \ldots, \alpha_m^{n-1}) = (a_m^0, \ldots, a_m^{n-1}), m = i, j)$$
$$(43) \qquad = \int_D P_{3n-1}^{\omega^{3n-2}}(D_1) \, dP^{3n-2}(\omega^{3n-2})$$
$$= \int_D \bar{P}^{h^n_{\omega^{3n-2}}}(D_2) \cdot \bar{P}^{h^n_{\omega^{3n-2}}}(D_3) \, dP^{3n-2}(\omega^{3n-2})$$
$$= P^{3n-2}(D)P(g_i^n = c | h_i^n = l)P(g_j^n = d | h_j^n = r),$$

where

$$D_1 = \{\omega_{3n-1} \in \Omega_{3n-1} : g^n(i, \omega^{3n-1}) = c, g^n(j, \omega^{3n-1}) = d\},$$
$$D_2 = \{\bar{\omega} \in \bar{\Omega} : \bar{g}^{h^n_{\omega^{3n-2}}}(i, \bar{\omega}) = c\},$$
$$D_3 = \{\bar{\omega} \in \bar{\Omega} : \bar{g}^{h^n_{\omega^{3n-2}}}(j, \bar{\omega}) = d\},$$

which means that the Markov conditional independence condition as formulated in (6) for random partial matching is satisfied.

Finally, we consider the step of random type changing for matched agents at time $n$. For $k \in S$, let $\nu_{kJ}$ be the Dirac measure at $k$ on $S$, that is, $\nu_{kJ}(r) = \delta_k^r$ for each $r \in S$. If $h^n(i, \omega^{3n-2}) = k$ for $k \in S$ and $g_i^n(\omega^{3n-1}) = c$ for $c \in S \cup \{J\}$, then the measure $\tau_i^{\omega^{3n-1}}$ in the definition of $Q_{3n}^{\omega^{3n-1}}$ on $S^I$ is simply $\nu_{kc}$.

For each agent $i \in I$, and for any $r, k \in S$, and $c \in S \cup \{J\}$, we have

$$P(\alpha_i^n = r, h_i^n = k, g_i^n = c)$$
$$= P^{3n}(\{(\omega^{3n-1}, \omega_{3n}) \in \Omega^{3n} : h_i^n(\omega^{3n-2}) = k,$$
(44)
$$g_i^n(\omega^{3n-1}) = c, \omega_{3n}(i) = r\})$$
$$= \int_{\{\omega^{3n-1} \in \Omega^{3n-1} : h_i^n(\omega^{3n-2})=k, g_i^n(\omega^{3n-1})=c\}} \nu_{kc}(r) \, dP^{3n-1}(\omega^{3n-1})$$
$$= \nu_{kc}(r) P(h_i^n = k, g_i^n = c),$$

which implies that (4) is satisfied.

Fix $i \neq j \in I$. Take any $(a_i^0, \ldots, a_i^{n-1}), (a_j^0, \ldots, a_j^{n-1}) \in S^n$, $k, l, r, t \in S$, and $c, d \in S \cup \{J\}$. Let $E$ be the set of all $\omega^{3n-1} \in \Omega^{3n-1}$ such that $g_i^n(\omega^{3n-1}) = c$, $g_j^n(\omega^{3n-1}) = d$, $h_i^n(\omega^{3n-2}) = k$, $h_j^n(\omega^{3n-2}) = l$, $(\alpha_i^0, \ldots, \alpha_i^{n-1})(\omega^{3n-3}) = (a_i^0, \ldots, a_i^{n-1})$ and $(\alpha_j^0, \ldots, \alpha_j^{n-1})(\omega^{3n-3}) = (a_j^0, \ldots, a_j^{n-1})$. Then, letting

$$E_1 = \{\alpha_i^n = r, \alpha_j^n = t, g_i^n = c, g_j^n = d, h_i^n = k, h_j^n = l,$$
$$(\alpha_m^0, \ldots, \alpha_m^{n-1}) = (a_m^0, \ldots, a_m^{n-1}), m = i, j\},$$

(44) implies that

$$P(E_1) = \int_E P_{3n}^{\omega^{3n-1}}(\{\omega_{3n} \in \Omega_{3n} : \omega_{3n}(i) = r, \omega_{3n}(j) = t\}) \, dP^{3n-1}(\omega^{3n-1})$$
(45)
$$= \int_E \nu_{kc}(r) \nu_{ld}(t) \, dP^{3n-1}(\omega^{3n-1})$$
$$= P(E) P(\alpha_i^n = r | h_i^n = k, g_i^n = c) P(\alpha_j^n = t | h_j^n = l, g_j^n = d),$$

which means that the Markov conditional independence condition as formulated in (7) for match-induced random type changing is satisfied.

Therefore, we have shown that $\mathbb{D}$ is a dynamical system with random mutation, partial matching and type changing that is Markov conditionally independent in types.

**Acknowledgments.** This work was initiated in July 2000 while the second author visited Stanford University. This work was presented in the annual meeting of the American Mathematical Society in Phoenix, Arizona on January 7–10, 2004 and in the "Workshop on Markov Processes and Related Topics" in Beijing Normal University on August 10–14, 2004. We are grateful to Bob Anderson, Ward Henson, Jerry Keisler, Horst Osswald, Peter Loeb and Jiang-Lun Wu for useful conversations, and to two anonymous referees of this journal for helpful comments.

GRADUATE SCHOOL OF BUSINESS  
STANFORD UNIVERSITY  
STANFORD, CALIFORNIA 94305  
USA  
E-MAIL: duffie@stanford.edu

DEPARTMENTS OF ECONOMICS AND MATHEMATICS  
NATIONAL UNIVERSITY OF SINGAPORE  
1 ARTS LINK  
SINGAPORE 117570  
E-MAIL: matsuny@nus.edu.sg